\documentclass[a4paper,12pt]{amsart}

\usepackage[usenames,dvipsnames]{xcolor}
\definecolor{brightcerulean}{rgb}{0.11, 0.67, 0.84}
\usepackage[pagebackref=true, colorlinks=false, allcolors=brightcerulean, allbordercolors=brightcerulean,
pdfauthor={J. Balakrishnan, N. Dogra, J.S. M\"uller, J. Tuitman, J. Vonk},
pdftitle={Quadratic Chabauty for Modular curves: Algorithms and examples}
]{hyperref}

\usepackage{url,bbold,dsfont}

\usepackage{tikz,enumerate,caption,stmaryrd,amsfonts,amssymb,systeme,comment,graphicx,colonequals,faktor,xfrac}
\usetikzlibrary{matrix,positioning,arrows,shapes,chains,calc,automata}
\usetikzlibrary{decorations.pathreplacing}
\usepackage[shortlabels]{enumitem}
\usepackage[T1]{fontenc}

\usepackage{comment}
\usepackage{graphicx}
\usepackage{epsfig}
\usepackage{fancyhdr}
\usepackage{mathrsfs}
\usepackage{libertine}

\renewcommand*{\backref}[1]{}
\renewcommand*{\backrefalt}[4]{%
  \ifcase #1 %
    \relax
  \or
    $\uparrow$#2.%
  \else
   $\uparrow$#2.%
  \fi%
}

\newcommand\cyr{%
  \renewcommand\rmdefault{cmr}%
  \renewcommand\sfdefault{wncyss}%
  \renewcommand\encodingdefault{OT2}%
  \normalfont\selectfont}
\DeclareTextFontCommand{\textcyr}{\cyr}

\definecolor{red}{rgb}{0.9,0,0}
\definecolor{purple}{rgb}{0.8,0,0.6}

\numberwithin{equation}{section}

\DeclareFontFamily{U}{wncy}{}
\DeclareFontShape{U}{wncy}{m}{n}{<->wncyr10}{}
\DeclareSymbolFont{mcy}{U}{wncy}{m}{n}
\DeclareMathSymbol{\Sha}{\mathord}{mcy}{"58}
\DeclareMathOperator{\NS}{\mathrm{NS}}

\DeclareMathOperator{\Q}{\mathbf{Q}}
\DeclareMathOperator{\QQ}{\mathbf{Q}}
\DeclareMathOperator{\F}{\mathbf{F}}
\DeclareMathOperator{\FF}{\mathbf{F}}
\DeclareMathOperator{\h}{H}
\DeclareMathOperator{\Z}{\mathbf{Z}}
\DeclareMathOperator{\ZZ}{\mathbf{Z}}
\DeclareMathOperator{\C}{\mathbf{C}}

\DeclareMathOperator{\PP}{\mathbf{P}}
\DeclareMathOperator{\HH}{\mathrm{H}}

\DeclareMathOperator{\ord}{\mathrm{ord}}
\DeclareMathOperator{\rk}{\mathrm{rk}}
\DeclareMathOperator{\SL}{SL}

\DeclareMathOperator{\et}{\scriptstyle \mathrm{\acute{e}t}}

\DeclareMathOperator{\dR}{\scriptstyle \mathrm{dR}}

\DeclareMathOperator{\fil}{\scriptstyle \mathrm{ Fil}}
\DeclareMathOperator{\cris}{\scriptstyle \mathrm{ cris}}
\DeclareMathOperator{\ns}{\scriptstyle \mathrm{ ns}}

\DeclareMathOperator{\rig}{\scriptstyle \mathrm{{rig}}}
\DeclareMathOperator{\st}{\scriptstyle \mathrm{st}}

\DeclareMathOperator{\Ext}{\mathrm{Ext}}

\DeclareMathOperator{\Hom}{\mathrm{Hom}}
\DeclareMathOperator{\End}{\mathrm{End}}

\DeclareMathOperator{\Gal}{\mathrm{Gal}}
\DeclareMathOperator{\Res}{\mathrm{Res}}

\DeclareMathOperator{\Fil}{\mathrm{Fil}}
\DeclareMathOperator{\Div}{\mathrm{Div}}
\DeclareMathOperator{\Tr}{\mathrm{Tr}}

\DeclareMathOperator{\GL}{\mathrm{GL}}

\DeclareMathOperator{\spl}{\mathrm{spl}}

\DeclareMathOperator{\lra}{\longrightarrow}

\newcommand{\Ker}{\mathrm{Ker}}

\renewcommand{\div}{\operatorname{div}}

\newcommand{\red}{\mathrm{red}}
\newcommand{\reg}{\mathrm{reg}}

\newtheorem{theorem}{Theorem}[section]
\newtheorem{prop}[theorem]{Proposition}

\newtheorem{lemma}[theorem]{Lemma}

\newtheorem{Definition}[theorem]{Definition}

\theoremstyle{remark}
\newtheorem{Remark}[theorem]{Remark}
\newtheorem{Example}[theorem]{Example}
\newtheorem{assump}[theorem]{Assumption}
\newtheorem{algorithm}[theorem]{Algorithm}

\usepackage{xypic}
\begin{document}
\dedicatory{Dedicated to the memory of Bas Edixhoven (1962 -- 2022)}

\title[Quadratic Chabauty: Algorithms and Examples]{ Quadratic Chabauty for Modular curves: \\ Algorithms and examples}

\author{Jennifer S. Balakrishnan}
\address{  Department of Mathematics  \&  Statistics,
Boston University, 665 Commonwealth Avenue, Boston, MA 02215, USA}
\author{ Netan Dogra}
\address{Department of Mathematics, King's College London, Strand, London, WC2R 2LS, UK  }
\author{J. Steffen M\"uller}
\address{ Bernoulli Institute, University of Groningen, Nijenborgh 9, 9747
AG Groningen, The Netherlands}
\author{Jan Tuitman} 
\address{ }
\author{Jan Vonk}
\address{Mathematical Institute, Leiden University, Niels Bohrweg 1, 2333 CA Leiden, The Netherlands} 
 
\thanks{JB was supported by NSF grant DMS-1945452,  the  Clare  Boothe
Luce  Professorship  (Henry  Luce  Foundation),  Simons Foundation grant
\#550023,  and a Sloan Research Fellowship.  ND was supported by a Royal
Society University Research Fellowship.  SM was supported by DFG grant MU
4110/1-1 and by NWO Grant VI.Vidi.192.106. JV was supported by ERC-COG Grant 724638
`GALOP' and Francis Brown, the Carolyn and Franco Gianturco Fellowship at
Linacre College (Oxford), and NSF Grant No. DMS-1638352, and NWO Grant
VI.Vidi.213.084 during various stages of this project.}

%=========================================================================
\begin{abstract} \setlength{\parskip}{1ex} \setlength{\parindent}{0mm}
We describe how the quadratic Chabauty method may be applied to determine the
  set of rational points on modular curves of genus $g>1$ whose Jacobians have
  Mordell--Weil rank $g$. This extends our previous work on the split Cartan curve of
  level 13 and allows us to consider modular curves that may have few known rational
  points or nontrivial local height contributions at primes of bad reduction.  We illustrate our algorithms with a
  number of examples where we determine the set of rational points on several modular
  curves of genus 2 and 3: this includes Atkin--Lehner quotients $X_0^+(N)$ of prime level
$N$, the curve $X_{S_4}(13)$, as well as a few other curves relevant to Mazur's Program B.
We also compute the set of rational points on the genus 6 non-split Cartan modular curve $X_{\ns} ^+ (17)$.\end{abstract}

\maketitle
%=========================================================================

%\tableofcontents

%=========================================================================
\section{Introduction}\label{sec:intro}
%=========================================================================
In this paper, we describe the current state of quadratic Chabauty--based algorithms for
the resolution of Diophantine equations arising from modular curves.  Here we consider the
usual modular curves associated to congruence subgroups of $\SL_2(\Z)$, as well as
Atkin--Lehner quotients thereof.

Recall the motivating question of the subject: let $E$ be an elliptic curve over a number field $K$. What are the possible ways for the Galois group $\Gal (\overline{K}/K)$ to act on the group of torsion points of $E$? Equivalently, what are the conjugacy classes of subgroups of $\GL _2 (\Z/N\Z)$ arising as images of the mod $N$ Galois representation $\rho _{E,N}$? 

By a theorem of Serre \cite{Ser72}, if $E$ is an elliptic curve without complex
multiplication, then for all primes $N\gg 0$, the representation $\rho _{E,N}$ is
surjective. Serre's uniformity question \cite{Ser72} asks whether this can be made
\textit{uniform} over $\Q $: is there an $N_0 $ such that, for all primes $N>N_0 $, if
$E/\Q $ is an elliptic curve without complex multiplication, then $\rho _{E,N}$ is
surjective? By a classification of maximal subgroups of $\GL _2 (\Z/N\Z)$, this amounts to
determining elliptic curves whose mod $N$ Galois representation is contained in a Borel subgroup,
the normaliser of a split Cartan subgroup, the normaliser of a non-split Cartan subgroup, or an `exceptional' subgroup (such that the projective image is $S_4 ,A_4,$ or $A_5$).

Mazur's Program B \cite{maz77} asks for all of the possible Galois actions on torsion
subgroups of elliptic curves without complex multiplication. This question includes
Serre's uniformity question but is more general. From a Diophantine perspective, it
roughly amounts to determining the rational points on all modular curves. 

Rouse
and Zureick-Brown~\cite{RZB} settled this in the context of 2-primary torsion
 and very recently, with Sutherland~\cite{RSZB}, studied this in the context of $\ell$-primary torsion for other primes $\ell$.  For
each prime, this produces a finite number of curves, the determination of whose rational
points would resolve the $\ell$-primary part of Mazur's question. In~\S\ref{subsec:s4}
and~\S\ref{subsec:dzb} we compute the rational points on 
four modular curves $X_{S_4}(13), X_{\ns}^+(17), X_{11}$, and $X_{15}$ arising in Mazur's Program B. In
particular, we show the following:
\begin{theorem}\label{T:XS4}
  We have $\#X_{S_4}(13)(\Q) = 4$. One of these points is a CM
  point, corresponding to discriminant $D=-3$. The other three are exceptional,
  with corresponding $j$-invariants listed in~\S\ref{subsec:s4}.
\end{theorem}
Here we call a non-cuspidal rational point~\textit{exceptional} if it does
not correspond to an elliptic curve with complex multiplication.
The curve $X_{S_4}(13)$ has genus~3. This completes the classification of elliptic curves $E/\Q$ and prime level $N>0$ such that
$\rho _{E,N}$ is contained in an exceptional subgroup.

We also determine the rational points on $X_{\ns}^+(17)$, the non-split Cartan modular curve of level 17, which is a genus 6 curve:
\begin{theorem}\label{T:Xns17}
  We have $\#X_{\ns}^+(17)(\Q) = 7$ and all of these points are CM
  , corresponding to discriminants $-3, -7, -11,$ $-12, -27, -28, -163$.
\end{theorem}
Theorems~\ref{T:XS4} and~\ref{T:Xns17} complete the classification of the possible 13-adic
and 17-adic images of Galois.

Moving beyond torsion points of elliptic curves over $\Q $, another interesting problem in the Diophantine geometry of modular curves is the determination of the set of rational points on the Atkin--Lehner quotient $$X_0^+(N) \colonequals X_0(N)/\langle w_N \rangle$$ of the modular curve $X_0 (N)$. In \cite{Gal02}, Galbraith asks whether, for all primes $N\gg 0$, the only rational points on $X_0 ^+ (N)$ are cusps or CM points. From a moduli perspective, this amounts to finding quadratic $\Q $-curves that are $N$-isogenous to their conjugates.  Dogra and Le Fourn \cite{DLF19} proved that the quadratic Chabauty set $X_0^+(N) (\Q _p )_2 $ is finite whenever the genus of $X_0^+(N)$ is larger than one. Hence it is natural to ask whether the methods of this paper can be used to give an algorithm for computing $X_0^+(N) (\Q _p )_2 $ for any $N$. In fact, in the range of $N$ we consider, finiteness of $X_0 ^+ (N)(\Q _p)_2$ follows from a criterion appearing in earlier work of Siksek \cite{Sik17}. 
Our computations described in~\S\ref{subsec:x0n+} prove the following result.
\begin{theorem}\label{T:X0N+}
  The only prime values $N$ such that the curve $X_0^+(N)$ is of genus 2 or 3 and has an
  exceptional rational point  are $N=73, 103, 191$. In particular for prime $N$, there are no exceptional
  rational points on curves $X_0^+(N)$ of genus 3.
\end{theorem}

All rational points in Theorem~\ref{T:X0N+} had already been found
by Galbraith~\cite{Gal99}.
\begin{Remark}
These computations were recently extended significantly by Ad\v{z}aga, Arul, Beneish,
  Chen, Chidambaram, Keller, and Wen \cite{AABCCKW}. They use the quadratic
  Chabauty method described in
  this paper to determine the set of rational points on all curves $X_0 ^+
  (N)$ of genus 4,5 and 6 and prime level $N$.  Arul and
  M\"uller~\cite{AM} also
  compute the rational points on $X_0^+(125)$ using the same method.
  Ad\v{z}aga, Chidambaram, Keller, and Padurariu \cite{ACKP} use several techniques, including quadratic Chabauty, to determine the set of rational points on the hyperelliptic Atkin--Lehner star quotient curves $X_0^*(N)$. 
\end{Remark}

Going further, one may wonder what the potential applications of these
algorithms are to \textit{non-modular} curves. The main stumbling block in
attempting such a generalisation is our running assumption on the
Mordell--Weil rank and Picard number of the Jacobian (see~\S\ref{S:setup}).
Since a generic curve has Picard number one, it is not clear how often one should expect a genus $g$ curve with Mordell--Weil rank $g$ to satisfy the quadratic Chabauty hypothesis.
Nevertheless, there are other interesting curves where one \textit{would} expect to get
some mileage out of such algorithms. The most obvious examples are (Atkin--Lehner
quotients of) Shimura curves. In particular, determining the set of rational points on the
(infinitely many) curves $X^D /\langle w_D\rangle$, in the notation of Parent--Yafaev \cite{PY}, would resolve a conjecture of Clark \cite{clarkthesis} (Parent and Yafaev determine the rational points for an infinite family of Shimura curves whose Jacobian contains a rank zero isogeny factor).

\subsection*{Acknowledgements} We are deeply indebted to Bas Edixhoven for his numerous generous insights on this subject.  It is a pleasure to thank Noam Elkies, Barry Mazur, Jeremy Rouse, Andrew Sutherland, and David Zureick-Brown for suggesting several  modular curves of interest, as well as many helpful discussions, which provided the impetus for this work.  We are grateful for the contributions of Nikola Ad\u zaga, Vishal Arul, Lea Beneish, Alex Best, Francesca Bianchi, Mingjie Chen, Shiva Chidambaram, Timo Keller, Nicholas Triantafillou, and Boya Wen 
in finding a number of bugs in earlier versions of our code. We would also like to thank
Michael Stoll for kindly providing an implementation of the Mordell--Weil sieve, on which
ours is based and
Nils Bruin for sharing another approach to determining $C_{188}(\Q)$.
We are grateful to Francesca Bianchi, David
Holmes, Timo Keller and Michael Stoll for helpful comments on an earlier
version of this article.
We thank the referee for a helpful and entertaining report.

%=========================================================================
\section{Quadratic Chabauty: Theory}\label{sec:qc_theory}
We give a brief overview of the quadratic Chabauty method. A more complete exposition can
be found in \cite{BBBLMTV19}, and we refer the reader to \cite{BD18,BDMTV19} for more
precise details and proofs. Our description is in terms of Galois representations and
filtered $\phi $-modules, but we note that recently, Edixhoven and Lido \cite{EL19} gave a
geometric version of quadratic Chabauty, which they used to determine  the  set of
rational points on the bielliptic modular curve $X_0(129)/\langle w_3, w_{43} \rangle$ of
genus~2. Duque-Rosero, Hashimoto, and Spelier~\cite{DRHS}
have related this approach to the one presented here and used this to  give
algorithms for geometric quadratic Chabauty for hyperelliptic curves.  Besser, M\"uller, and Srinivasan~\cite{BMS21} have also given an
alternative approach to the quadratic Chabauty method based on a new
construction of $p$-adic heights on abelian varieties via $p$-adic Arakelov
theory.

\par An early version of the method appeared in work of Kim \cite{Kim10b, BKK10}, where Massey products were used to construct a locally analytic function, vanishing on the set of integral points of an elliptic curve of rank $1$. These functions were interpreted as height functions, extending the method, in Balakrishnan--Besser \cite{BB15} and Balakrishnan--Besser--M\"uller \cite{BBM16}. It was extended to its current form in Balakrishnan--Dogra \cite{BD18}, where a systematic use of Nekov\'a\v{r}'s theory of $p$-adic heights suggested a streamlined approach towards a very general class of curves allowing an abundance of geometric correspondences. It was carried out to determine the set of rational points on $X_{\rm s}^+(13)$, the split Cartan curve of level 13, in \cite{BDMTV19}.

\par \textbf{Remark. }This method fits into the vastly more general framework developed by Kim \cite{Kim05,Kim09}, elaborating on the idea of studying rational points on curves through path torsors of the \'etale fundamental group, suggested by Grothendieck's section conjecture. The approach discussed here represents an effective way to make this theory computable and applicable to a variety of examples. It is, however, important to note that different quotients of the fundamental group have been successfully used for this purpose, see for instance \cite{BD19}. Finally, although we restrict our attention to the base field $\Q$, suitable versions exist over number fields, see \cite{BD18, BD19, BBBM}.

\subsection{Rational points and global heights. }\label{S:setup}
Consider a smooth projective curve $X_{\Q}$ of genus $g \geq 2$ whose Jacobian $J$ has rank $r=g$. 
We also assume that the abelian logarithm induces an isomorphism 
\begin{equation}\label{eq:log_isom}
\log\colon J(\Q)\otimes\Q_p \to \HH^0(X_{\Q_p},\Omega^1)^{\vee}
\end{equation}
and that $X(\Q)$ is non-empty, so we may
choose a base point $b$ in $X(\Q)$. Suppose that the N\'eron--Severi rank $\mathrm{rk}_{\Z} \mathrm{NS}(J)$ is at least $2$, so that there exists a nontrivial class 
\[
\label{eqn:Z-choice}
Z \ \in \ \mathrm{Ker}\left(\mathrm{NS}(J) \longrightarrow \mathrm{NS}(X) \simeq \Z \right).
\]
As explained in Balakrishnan--Dogra \cite[Lemma 3.2]{BD18}, we can attach to any such choice of $Z$
a suitable quotient $U_Z$ of the $\Q_p$-pro-unipotent fundamental group
of $X_{\bar{\Q}}$, which via a twisting construction by path torsors, gives rise to a certain family of Galois representations
\[
\begin{array}{ccc}
X(K) & \longrightarrow & \left\{ G_K \to \mathrm{GL}_{2g+2}(\Q_p) \right\}/\sim \\
x & \longmapsto & \mathrm{A}(x) \colonequals \mathrm{A}_Z(b,x)
\end{array}
\]
where $K \in \{ \Q, \Q_p \}$ and $G_K$ is the absolute Galois group of
$K$. We refer the reader to \cite[\S 5.1]{BD18} 
for the details of this construction (in particular for the equivalence relation), and merely recall here that with respect to a suitable choice of basis, the representation $\mathrm{A}(x)$ is lower triangular, of the form
\begin{equation}
\label{eqn:matrix_Ex}
g \in G_K \ \longmapsto \ 
\left( 
\begin{matrix}
1 & & \\
\alpha(g) & \rho_V(g) & \\
\gamma(g) & \beta(g) & \chi_p(g) 
\end{matrix}
\right)
\end{equation}
where 
\[
\rho_V \colon G_K \longrightarrow \mathrm{GL}_{2g}(\Q_p)
\]
is a frame for the Galois action on the $p$-adic \'etale homology $V =
\HH^1_{\et}(X_{\overline{K}}, \Q_p)^{\vee}$, and $\chi_p\colon G_K \to \Q_p^{\times}$ is the $p$-adic cyclotomic character. Representations of this form, which admit a $G_K$-stable filtration with graded pieces $\Q_p(1), V, \Q_p$, are referred to as \textit{mixed extensions}, see \cite[\S 3.1]{BDMTV19}.

\par The theory of $p$-adic heights due to Nekov\'a\v{r} \cite[\S 2]{Nek93} attaches to any mixed extension $M$ a $p$-adic height $h(M)$. When applied to the family of mixed extensions $\mathrm{A}(x)$,  this results in a map 
\[
\label{eqn:global-height}
h \colon X(\Q) \longrightarrow \Q_p.
\]
The algebraic properties of this map lie at the heart of the quadratic Chabauty method. Most notably, the method relies on the following two facts:
\begin{itemize}
\item The $p$-adic height is a bilinear function of the pair of cohomology classes $([\alpha ], [\beta ])$ associated to the vectors appearing in \eqref{eqn:matrix_Ex}.
\item It
  decomposes as a sum of local height functions $h_v$ defined locally at every finite place $v$. 
\end{itemize}

\subsection{Local decomposition. }
\label{subsec:local_decomp}
We now discuss in more detail the decomposition of the global $p$-adic height $h$ described above, as a sum of local height functions
\[
h_v \colon X(\Q_v) \ \longrightarrow \ \Q_p. 
\]
The nature of these local height functions is as follows:
\begin{enumerate}
\item \textbf{The case $v \neq p$:} It follows from Kim--Tamagawa
  \cite[Corollary 0.2]{KT08} 
that the function $h_v$ has finite image, in the sense that there exists a finite set $\Upsilon_v$ such that
\[
h_v \colon X(\Q_v) \ \longrightarrow \ \Upsilon_v \subset \Q_p.
\]
\item \textbf{The case $v=p$:} The map $h_p$ is locally analytic and has a simple description in terms of linear algebra data of the filtered $\phi$-module 
\[
M(x) \colonequals \left( \mathrm{A}(x) \otimes_{\Q_p} \mathrm{B}_{\cris}\right)^{G_{\Q_p}},
\]
where $\mathrm{B}_{\cris}$ is Fontaine's crystalline period ring. A crucial
    point in the method of quadratic Chabauty is that the definition of the family of Galois representations $\mathrm{A}(x)$ comes from a motivic quotient of the fundamental group of $X$, and non-abelian $p$-adic Hodge theory yields an analogous de Rham realisation in the form of a filtered connection $(\mathscr{M},\nabla)$ on $X$ with a Frobenius structure, together with an isomorphism of filtered $\phi$-modules
\[
x^* \mathscr{M} \simeq M(x)
\]
(see \cite[\S 5]{BDMTV19}). We have a pair of elements $\pi _1 (M(x))$ and $\pi _2 (M(x))^\vee (1)$ of $\HH^0 (X_{\Q _p },\Omega )^\vee $ associated to the filtered $\phi $-module $M(x)$, via the isomorphism 
\[
\Ext ^1 _{\fil ,\phi }(\Q _p  ,\HH ^1 _{\dR}(X_{\Q _p })^\vee )\simeq \HH ^0 (X_{\Q _p },\Omega )^\vee .
\]
\end{enumerate}
\subsection{Finiteness. }
The decomposition $h = \sum_v h_v$ can be used to leverage the bilinear nature of $h$
against the properties of the functions $h_v$. By (1) in~\S2.2, we know that there exists a finite set
$\Upsilon = \Upsilon_Z \subset \Q_p$ such that 
\begin{equation}
\label{eqn:QCpair}
h(x) - h_p(x)  \ \in \ \Upsilon
\end{equation}
for any $x$ in $X(\Q)$. 
In Section~\ref{sec:qc_algorithms}, we describe how the terms in this equation may be
computed explicitly.
\begin{itemize}
\item The set $\Upsilon $ is given by $\{\sum _v  \epsilon _v :\epsilon _v \in \Upsilon _v \}$, where the sum is over primes of bad reduction, and $\Upsilon _v $ is the set of values of $h_v (x)$ for $x\in X(\Q _v ).$
  For $v\ne p$, the map $h_v$ is made more explicit in
    \S\ref{subsec:betts_dogra} 

  using the results of Betts--Dogra \cite{BeD19} to compute $\Upsilon_v$ when a regular semi-stable model $\mathcal{X}$ is known. The map $h_v$ factors through the reduction map to the irreducible components of the special fibre of $\mathcal{X}$. 
\item The map $h_p$ may be computed using  \cite[\S\S4,5]{BDMTV19}, where it is explained how the universal properties of the bundle $\mathscr{M}$ rigidify the (known) structures on the graded pieces, enough to allow us to compute them explicitly, see \S\ref{subsec:height_at_p}. 
\item Using the isomorphism \eqref{eq:log_isom}, we may view the global height as a pairing
\[
h \colon \HH^0(X_{\Q_p},\Omega^1)^{\vee} \otimes \HH^0(X_{\Q_p},\Omega^1)^{\vee} \ \longrightarrow \ \Q_p.
\]
Using global information, such as an abundance of global points $x \in X(\Q)$ if
    available, we can solve for the height pairing. This is discussed in \S\ref{subsec:global_height}, where we  also explain what to do when too few rational points are available. 
\end{itemize}

\par Via the above, the map $h$ may be extended to a bilinear map
\begin{equation}\label{eq:hdef}
h\colon X(\Q_p) \to \Q_p\,;\qquad x\mapsto h(\pi _1 (\mathrm{A}(x)),\pi _2
(\mathrm{A}(x))^\vee (1))\,.
\end{equation}
The resulting map
\begin{equation}\label{eq:rho}
\rho = h - h_p \colon X(\Q_p) \ \longrightarrow \ \Q_p
\end{equation}
is known to be Zariski dense on every residue disk.
We call $\rho$ a {\em quadratic Chabauty function}, and we write
$\rho_Z$ if we want to emphasise the dependence on $Z$.
Hence \eqref{eqn:QCpair} implies that $X(\Q)$ is finite. Moreover, the computable nature of the quantities involved in \eqref{eqn:QCpair}, discussed at length in the next section, allows us to explicitly determine a $p$-adic approximation of the finite set
\[
\{ x\in X(\Q _p )\colon h(x)-h_p (x)\in \Upsilon \} \supset X(\Q )\,.
\] 
As explained in \cite[Proposition 5.5]{BD18}, this finite set contains the
Chabauty--Kim set $X(\Q _p )_2$. In particular, a proof that this set
equals $X(\Q )$ gives a verification of Kim's conjecture \cite[Conjecture
3.1]{BDCKW} for the curve $X$ (we refer the reader to \cite[Definition
2.7]{BDCKW}
for the definitions of the set $X(\Q _p )_2 $).

%=========================================================================
\section{Quadratic Chabauty: Algorithms}\label{sec:qc_algorithms}
%=========================================================================

In this section, we discuss the computation of the three ingredients outlined above: 
\begin{enumerate}
\item The local height function $h_v$ for $v$ away from $p$, which is described in \S\ref{subsec:betts_dogra} using the techniques in Betts--Dogra ~\cite{BeD19}, given a regular semi-stable model at $v$.
\item The height function $h_p$, whose computation using the techniques of~\cite{BDMTV19} is described in
  \S\ref{subsec:height_at_p} 
\item The determination of the global height pairing $h$, described in \S\ref{subsec:global_height} using rational divisors as input in the absence of a supply of rational points on the curve. 
\end{enumerate}

\par Our contribution in this paper lies mainly in (1) and (3), which reflect general
features of the method of quadratic Chabauty that were not needed for the curve $X_{\rm
s}^+(13)$ treated in ~\cite{BDMTV19}. In addition, we discuss some computational
techniques to further automate the method of quadratic Chabauty to work for a wide class
of modular curves. This includes the Mordell--Weil sieve, which is used to attempt to further refine the finite set of local points in the output to the true set of rational points $X(\Q)$. 

\begin{Remark}\label{R:choices}The global height depends on the choice (which we fix henceforth) of 
\begin{itemize}
  \item a nontrivial continuous id\`ele class character $\chi\colon \mathbf{A}_{\Q}^\times / \Q^{\times}
  \lra \Q_p$ ramified at $p$; 
\item a splitting $s\colon V_{\dR}/\Fil^0 V_{\dR} \lra V_{\dR}$ of the Hodge filtration, where
  $$V_{\dR} = \mathrm{D}_{\cris}(V) =  \HH^1_{\dR}(X_{\Q_p})^{\vee}\,.$$
\end{itemize}
We also fix differentials $\omega_{0},\ldots,\omega_{2g-1}$ of the second kind whose
classes form a symplectic basis of $\HH^1_{\dR}(X_{\Q_p})$ with respect to
the cup product, such that
$\omega_{0},\ldots,\omega_{g-1}$ generate $\HH^0(X_{\Q_p},\Omega^1)$.
\end{Remark}

%%%%%%%%%%%%%%%%%%%%%%%%%%%%%%%%%%%%%%%%%%%%%%%%%%%%%%%%%%%%%%%%%%%%%%%%
\subsection{Local heights away from $p$}\label{subsec:betts_dogra}
Let $\ell \neq p$ and let $F$ be an endomorphism
of $J$ whose class $Z$ lies in 
$\mathrm{Ker}\left(\mathrm{NS}(J) \to \mathrm{NS}(X) \right)$. 
In ~\cite{BeD19}, a description of the map
\[
h_\ell \colon X(\Q _{\ell }) \lra \HH^1(G_\ell, U_Z) \lra  \HH^1(G_\ell,
  \Q_p(1)) \lra     \Q _{p}
\]
associated to $F$ and $\chi $ is given, in terms of harmonic analysis on the reduction
graph in the sense of Zhang ~\cite{zhang}.

To explain the result, we introduce some notation. Over some finite
extension $K/\Q _{\ell }$, the curve $X$ admits a regular semistable model
$\mathcal{X}_{\mathrm{reg}}/\mathcal{O}_K$, and a stable model
$\mathcal{X}_{\mathrm{st}}/\mathcal{O}_K$. Let $\Gamma _{\mathrm{reg}}$ and
$\Gamma _{\mathrm{st}}$ denote the dual graphs of the special fibres of
these models. Recall that the \textit{dual graph} of the special fibre
is by definition the graph \footnote{Here we follow the convention that
graphs are allowed multiple edges between two vertices, and loops (i.e. an
edge whose endpoints are equal).} whose vertices are the irreducible
components of the special fibre, and whose edges are the singular points of
the special fibre. The endpoints of an edge $e$ are defined to be the
irreducible components containing the point (by semistability, a singular
point $e$ lies on at most two irreducible components). By regularity, we have a reduction map
\[
  \red\colon X(\Q _{\ell }) \lra V(\Gamma _{\mathrm{reg}})
\]
from $X(\Q _{\ell })$ to the vertices of the dual graph $\Gamma _{\mathrm{reg}}$. 

The definition is the natural one: given $x\in X(\Q _{\ell })$, there is a
unique extension to an $\mathcal{O}_K$-section $x\in
\mathcal{X}_{\mathrm{reg}}(\mathcal{O}_K )$. Let $k$ be the residue
field of $\mathcal{O}_K$. By regularity, the specialisation of $x$ to $k$ lies on a unique irreducible component of $\mathcal{X}_{\mathrm{reg},k}$.

We may give $\Gamma _{\mathrm{reg}}$ and $\Gamma _{\mathrm{st}}$ the structure of \textit{rationally metrised} graphs (i.e. graphs whose edges $e$ have associated lengths $\ell (e)\in \mathbb{\Q }_{>0}$) by defining the length of an edge $e$ to be $i(e)/r$, where $i$ is the intersection multiplicity of the corresponding singular point and $r$ is the ramification degree of $K/\Q _{\ell }$. 

Choose an orientation of the edges of $\Gamma\colonequals \Gamma _{\mathrm{st}}$, so that each $e\in E(\Gamma )$ has a source $s(e)$ and target $t(e)$ in $V(\Gamma )$. We define the (rational) homology of $\Gamma $, $H_1 (\Gamma )\subset \Q E(\Gamma )$, to be the kernel of the map
\[
s-t\colon \Q E(\Gamma )\to \Q V(\Gamma ),
\]
where $\Q E(\Gamma )$ and $\Q V(\Gamma )$ are the free $\Q $-vector spaces generated by $E(\Gamma )$ and $V(\Gamma )$ respectively.

Define $\Gamma _{\Q }$ to be the set of
points on $\Gamma$ whose distance from a vertex is rational: formally,
\[
  \Gamma _{\Q} =\sqcup _{e\in E(\Gamma _{\st})}\{e\}\times ([0,\ell (e)] \cap \Q ) /\sim \, ,
\]
where the equivalence relation is that $(e_1 ,1 )\sim (e_2 ,0)$ whenever $t(e_1 )=s(e_2
)$. Since $\mathcal{X}_{\reg}$ is obtained from $\mathcal{X}_{\st}$ by taking each singular point (corresponding to an edge $e$) and blowing up $i(e)$ times, we have an inclusion $V(\Gamma _{\reg})\subset \Gamma
_{\Q }$ (in the terminology of \cite[3.7.1]{BeD19}, we may view $\Gamma _{\reg}$ as a rational subdivision of $\Gamma _{\st }$). 
In this way we can think of the reduction map $\red
$ as a map from $X(K)$ to $\Gamma _{\Q }$,
see~\cite[Definition~1.3.1]{BeD19}. The rationally metrised graph we
obtain is independent of the choice of extension over which $X$ acquires
stable reduction \cite[Proposition 2.6]{CR91}, and in fact there is an
equivalent definition of $\Gamma _{\Q }$ as the limit of the dual graphs of
special fibres of regular semistable models of $X_L$ over all finite
extensions $L$ of $K$ (see \cite[\S 2]{CR93}).

In \cite[Lemma 12.1.1]{BeD19}, a map 
\[
j_{\Gamma }\colon\Gamma _{\Q }\to \Q _p
\]
is defined such that $h_\ell =c\cdot j_{\Gamma }\circ \red $, where $c$ is
a constant. 
The map $j_{\Gamma }$ is defined in terms of the Laplacian operator associated to $\Gamma _{\mathrm{st}}
$, which we now define. We say a function
\[
\Gamma _{\Q} \to \Q _p
\]
is \textit{piecewise polynomial} if on each edge it is the restriction of a polynomial
function $\Q \to \Q _p $. As in~\cite[Definition~7.2.2]{BeD19}, we define the \textit{Laplacian} $\nabla ^2 (g)$ of a piecewise
polynomial function $g\colon \Gamma _{\Q }\to \Q _p$ to be the formal sum
\[
-\sum _{e\in E(\Gamma )} g '' (x_e )\cdot e +\sum _{v\in V(\Gamma )}(\sum _{s(e)=v}g'(0)-\sum _{t(e)=v}g'(1))\cdot v.
\]
Here we write the function $g$ restricted to the edge $e$ as a polynomial in $\Q _p [x_e ]$
for notational simplicity, where $x_e$ is the inclusion from the edge
$e$, thought of as a line segment $[0,\ell (e)]\cap \Q $, into $\Q $. Hence we have  $$\nabla ^2 (g) \in \bigoplus _{e\in E(\Gamma )}\Q _p [x_e
]\cdot e \oplus \bigoplus _{v\in V(\Gamma )}\Q _p \cdot v \,.$$ 
The Laplacian is linear on piecewise polynomial functions, and its kernel consists of constant functions. Thus $g$ is uniquely determined by $\nabla ^2(g)$ and its value at one point.

In \cite{BeD19}, an explicit construction is given of a piecewise
polynomial function that corresponds, via $\red $, to the local height
function we wish to compute. Recall that $F$ is an element of $\End
(J)\otimes \Q _p$ whose image in $\NS (J)$ lies in the kernel of $\NS
(J)\to \NS (X)$, and $b\in X(\Q )$ is a rational point.
\begin{theorem}[{\cite[Theorem 1.1.2, Lemma 12.1.1,  and Corollary 12.1.3]{BeD19}}]\label{thm:BeD}
Let $\Gamma $ be the dual graph of $X$ corresponding to a regular semi-stable model of $X$
  over $\mathcal{O}_K$, where $K/\Q _{\ell }$ is a finite extension. Let
  $\red \colon X(\Q _{\ell } )\to V(\Gamma )$ be the reduction map. For an
  irreducible component $X_w$ 
  of the special fibre of the regular semistable model, let $V_p (X_w )$ denote the $\Q _p $-Tate module of its Jacobian.
The morphism $j_{\Gamma }$ is the unique piecewise polynomial function
\[
j_\Gamma \colon \Gamma _{\Q }\to \Q _p
\]
satisfying $j_{\Gamma }(\red (b))=0$ and $\nabla ^2 (j_{\Gamma })=\mu _F $, where
\[
\mu _F \colonequals \sum _{e\in E(\Gamma )}\frac{1}{\ell (e)} e^* F(\pi (e)) \cdot e+\frac{1}{2}\sum _{w\in V(\Gamma )}\Tr (F|V_p (X_w ))\cdot w .
\]
\end{theorem}
Here, the morphism $\pi $ is by definition the orthogonal projection
\[
\Q E(\Gamma )\to \HH_1 (\Gamma ,\Q )
\]
with respect to the pairing $e\cdot e'=\delta _{ee'}$ on $\Q E(\Gamma
)$, and $e^* $ is the functional $\Q E(\Gamma )\to \Q $ projecting
onto the $e$ component. 
Recall (e.g. \cite[12.3.7]{SGA7}) that $V_p (X)$ admits a $G_K$-stable filtration
\[
V_p (X)=W_0 V_p (X)\supset W_1 V_p (X)\supset W_2 V_p (X)\supset W_3 V_p (X)=0,
\]
and we have isomorphisms of $G_K$-representations
\begin{align*}
\mathrm{gr}^W _0 V_p (X) & \simeq H_1 (\Gamma )\otimes \Q _p , \\
\mathrm{gr}^W _1 V_p (X) & \simeq \oplus _{w\in V(\Gamma )}V_p (X_w ), \\
\mathrm{gr}^W _2 V_p (X) & \simeq H_1 (\Gamma )^* \otimes \Q _p (1).
\end{align*}
The action of $F$ on $V_p (X)$ preserves this filtration since
it is a morphism of Galois representations, and hence induces an action of
$F$ on the weight $-1$ part of $V_p (X)$, which is isomorphic to $\oplus _w V_p (X_w )$. Although the action of $F$ need not respect the direct sum decomposition, the decomposition 
\[
\End (\oplus _w V_p (X_w ))\simeq \oplus _{w_1 ,w_2 }\Hom (V_p (X_{w_1 }),V_p (X_{w_2} ))
\]
implies that we can define $\Tr (F|V_p (X_w ))$ as the trace of the $\End
(V_p (X_w ))$-component of $F$.

To determine the possible local heights, it suffices to compute the action of $F$ on $\HH_1 (\Gamma )$ and on $V_p (X_v )$. In this paper, we do not discuss methods for the algorithmic computation of the action of $F$ on $\HH_1 (\Gamma )$, but algorithms for these computations in the case when the curve $X$ is hyperelliptic will be discussed in forthcoming joint work of the first, second and fifth authors with David Corwin, Sachi Hashimoto, Benjamin Matschke, Oana Padurariu, Ciaran Schembri, and Tian Wang.

As we explain in Section \ref{subsec:nontriv_away}, one can sometimes use partial
information deduced from Theorem~\ref{thm:BeD} to determine the possible
local heights without computing the action of $F$ on $\HH_1 (\Gamma)$ (for
example, if one has enough rational points on $X$ that are suitably independent in $J(\Q)$ and $\Gamma _{\Q }$). 

\begin{Example}\label{E:h2_188}
  One example for which this strategy succeeds is
  the curve $C_{188}/\Q$ defined by the equation $y^2 = x^5 - x^4 + x^3 + x^2 - 2x + 1$, as described in Example \ref{E:188}.  This curve does not have semistable reduction over $\Q_2$. 
 Over $K=\Q_2[\sqrt[3]{2}]$, we find a regular semistable model $\mathcal{X}_{\mathrm{reg}}$ whose special fibre consists 
  of two genus~1 curves that do not intersect and a genus~0 curve intersecting both of them transversely in a
  unique point each.
  We did not manage to obtain this information using any of the existing software packages for computing regular or
  semistable models, such as {\tt Magma}'s
  {\tt RegularModel} or the {\tt SageMath} package {\tt MCLF}~\footnote{MCLF
  can be used to show that there is a semistable model with three components,
  two of genus~1 and one of genus~0. It also lists equations for their function
  fields, but this information does not suffice for our purposes.}) 
  Therefore we computed this model by hand, using a standard (but
  tedious) sequence of blow-ups.
  
  Hence the metric graph $\Gamma_{\reg}$
  is a line segment and the image of $C_{188}(\Q _2 )$ in $\Gamma _{\Q }$ consists of three points on this line. The two edges of $\Gamma_{\reg}$ both have length $1/3$.
  In this case, since $\Gamma$
  has trivial homology, the function $j_{\Gamma }$ is affine linear, so it is uniquely determined by evaluating it at two distinct points.
  We use this to compute
  the rational points on $C_{188}$ in Example \ref{E:188}.
\end{Example}

%%%%%%%%%%%%%%%%%%%%%%%%%%%%%%%%%%%%%%%%%%%%%%%%%%%%%%%%%%%%%%%%%%%%%%
\subsection{Local heights at $p$}\label{subsec:height_at_p}
We discuss the local height component
\[
h_p \colon X(\Q_p) \ \lra \ \Q_p\,,
\]
which appeared in \cite[\S 5]{BDMTV19}. Recall that $h_p$ is a locally analytic function, described in terms of the filtered $\phi$-module $M(x)$ discussed in \S\ref{subsec:local_decomp}. Concretely, we may find two unipotent isomorphisms 
\[
\lambda^{\star}(x) \colon \ \Q_p \oplus V_{\dR} \oplus \Q_p(1)  \ \stackrel{\sim }{\lra } \ M(x), \quad\qquad \mbox{for} \ \ \star \in \{ \phi, \Fil \}
\]
where $\lambda^{\phi}$ respects the Frobenius action and $\lambda^{\Fil}$ respects the Hodge filtration, which with respect to a suitable basis for $M(x)$ may be represented in ($1+2g+1$)-block matrix form as
\begin{equation}\label{eqn:definitions-abcs}
\lambda^\phi(x) =\left(
\begin{array}{ccc} 
1 & 0 & 0 \\
\boldsymbol{\alpha }_{\phi} & 1 & 0 \\
\gamma_{\phi} & \boldsymbol{\beta }^{\intercal}_{\phi} & 1 \\
\end{array}
\right), \qquad 
\lambda^{\Fil}(x) = \left( 
\begin{array}{ccc}
1 & 0 & 0 \\
\boldsymbol{\alpha }_{\Fil}  & 1 & 0 \\
\gamma_{\Fil} & \boldsymbol{\beta }^{\intercal}_{\Fil} & 1 \\
\end{array}
\right)
\end{equation}
(see \cite[\S 5.3]{BDMTV19} and \cite[\S 4.5]{BDMTV19} respectively).
The isomorphism $\lambda^\phi$ is uniquely determined, whereas $\lambda^{\Fil}$ is only well-defined up to the stabiliser of the Hodge filtration $\Fil^0$.  A suitable choice gives  $\boldsymbol{\alpha }_{\Fil} =0$. 

\par 
The splitting $s$ of the Hodge filtration (see
Remark~\ref{R:choices}) defines idempotents $s_1,s_2$ on $V_{\dR}$ with images $s(V_{\dR}/\Fil^0 V_{\dR})$ and $\Fil^0 V_{\dR}$ respectively, with respect to which the local height at $p$ is
\begin{equation}\label{hpformula}
h_p (x)=
 \gamma_{\phi}
 - \gamma_{\Fil}
 -\boldsymbol{\beta}^\intercal_{\phi} \cdot s_1 (\boldsymbol{\alpha}_{\phi})
 - \boldsymbol{\beta }^\intercal_{\Fil} \cdot s_2 (\boldsymbol{\alpha}_{\phi})
\end{equation}
by~\cite[Equation~(17)]{BDMTV19}.

\par In \cite{BDMTV19} we outline a method to compute these quantities
explicitly as functions of the local point $x$ in $X(\Q_p)$, which exploits
the existence of the connection $(\mathscr{M},\nabla)$ discussed in
\S\ref{subsec:local_decomp}. The Hodge filtration and Frobenius structures
of this bundle are characterised by suitable universal properties,
discussed at length in \cite[\S\S 4--5]{BDMTV19}. We have made the
algorithms for the computation of $h_p$ more general and streamlined
and have added a precision analysis in Section~\ref{sec:prec} but did not make further contributions to this part of the method beyond what is already contained in \textit{loc. cit.}

%%%%%%%%%%%%%%%%%%%%%%%%%%%%%%%%%%%%%%%%%%%%%%%%%%%%%%%%%%%%%%%%%%%%%%%%
\subsection{The global height pairing}\label{subsec:global_height}
One key step in the construction of a quadratic Chabauty function is to write the global height pairing $h$ in terms of a basis of the space
of bilinear pairings on $\HH^0(X_{\Q_p}, \Omega^1)^{\vee}$. 
In \cite{BDMTV19}, we had as a working hypothesis that our curve $X$ had
sufficiently many rational points, in the following sense:
For $x\in X(\Q_p)$, the Galois representation $\mathrm{A}(x)$ can be projected onto $\HH^1_f(G_T, V)$
(respectively $\HH^1_f(G_T, V^*(1))$), where $G_T$ is the maximal
quotient of $G_{\Q}$ unramified outside $T =\{p\}\cup\{$bad primes for
$X\}$.  With respect to the dual basis
$\omega^\ast_0,\ldots,\omega^\ast_{g-1},$ the image is the vector $\alpha$ (respectively
$\beta)$ in~\eqref{eqn:matrix_Ex}. Both of these cohomology groups are isomorphic, under
our running assumptions, to $\HH^0(X_{\Q_p}, \Omega^1)^{\vee}$, so we obtain 
$$\pi(\mathrm{A}(x))= \left(\pi_1(\mathrm{A}(x)), \pi_2(\mathrm{A}(x))\right)  \in \HH^0(X_{\Q_p}, \Omega^1)^{\vee}\times \HH^0(X_{\Q_p}, \Omega^1)^{\vee}\,.$$
Suppose that we can find a basis of $\HH^0(X_{\Q_p}, \Omega^1)^{\vee} \otimes
\HH^0(X_{\Q_p}, \Omega^1)^{\vee}$ consisting of elements of the form $\pi(\mathrm{A}_Z(b,x))$,
where the $Z$ are cycles on $J$ pulling back to degree~0 cycles on $X$, and
the $x$ are rational points on $X$. Then we can compute the coefficients of $h$ in terms of the dual basis by evaluating $h_p(\mathrm{A}_Z(b,x))$ (and, if necessary, $h_{\ell}(\mathrm{A}_Z(b,x))$ for primes $\ell\ne p$). With this choice of basis, the extension of $h$ to a locally analytic function $h\colon X(\Q_p)\to\Q_p$ is immediate.

The number of required rational points can be reduced by working with symmetric heights
that are $\End(J)$-equivariant.
By
the latter we mean that $h(f(x), y) = h(x, f(y))$ for all $f\in \End(J)$, using~\eqref{eq:log_isom}. 
This holds if the splitting $s$ of the Hodge filtration on $V_{\dR}$ commutes with
$\End(J)$ and has the property that $\ker(s)$ is isotropic with respect to the cup product
(see~\cite[\S 4.11]{Nek93} and~\cite[\S 4.1]{BD19}). 
For instance, if $p$ is a prime of ordinary reduction for the Jacobian, then the height associated to the
unit root splitting (see Remark~\ref{rem:unitroot}) is symmetric and $\End(J)$-equivariant.
Henceforth we shall assume that 
$s$ satisfies these assumptions, and we say that $X$ has {\em sufficiently many rational
points} if the approach outlined above succeeds.

\subsubsection{Heights on the Jacobian}\label{subsec:coleman-gross} 
If our curve does not have sufficiently many rational points in the above sense, then, in
light of~\eqref{eq:log_isom}, it
is natural to solve for the height pairing using rational points on the Jacobian. In this case, we do not have an
algorithm at our disposal to compute $h$ using Nekov\'a\v r's construction, but we can use the equivalence between this
construction and that of Coleman and Gross~\cite{CG89}, proved by
Besser~\cite{Bes04}. In the case when the curve is hyperelliptic and
given by an odd degree model over $\Q_p$ (but see Remark~\ref{R:stevan}), we can further use the algorithm of Balakrishnan--Besser \cite{bb:heights, bb:heightserrata}. In the discussion that follows, we will assume that we are in this situation.
We will also assume that we know $g$ independent points on the Jacobian.

Recall from Remark~\ref{R:choices} that we have fixed a a continuous id\`ele class character $\chi\colon \mathbf{A}_{\Q}^\times / \Q^{\times}
  \lra \Q_p$ ramified at $p$ and 
a splitting $s\colon V_{\dR}/\Fil^0 V_{\dR} \lra V_{\dR}$ of the Hodge
filtration on $V_{\dR} = \HH^1_{\dR}(X_{\Q_p})^{\vee}\,.$ The latter
corresponds to a subspace $W\subset H^1_{\dR}(X_{\Q_p})$, complementary to
the image of $\HH^0(X_{\Q_p},\Omega^1)$. 
With respect to these choices, 
Coleman and Gross define the {\em local $p$-adic height pairing} $h_v(D_1,D_2) \in \Q_p$ at a finite prime $v$
 for divisors $D_1,D_2\in \Div^0(X_{\Q_v})$ with disjoint support.
The local pairing is
bi-additive, and we have $h_v(D_1,D_2) = \chi_v(f(D_2))$ if $D_2 = \div(f)$ is
principal. For $v\ne p$, the pairing $h_v$ is also symmetric; $h_p$ is
symmetric if and only if $W$ is isotropic with respect to the cup product
pairing, which we will assume from now on. Moreover, for $D_1,D_2\in \Div^0(X)$ with disjoint support, only finitely many
$h_v(D_1,D_2)\colonequals h_v(D_1\otimes \Q_v ,D_2\otimes \Q_v)$ are nonzero. Therefore $h\colonequals \sum_v h_v$ defines 
a symmetric bilinear pairing $h\colon J(\Q)\times J(\Q)\to \Q_p$ (see \cite[\S 6]{CG89}).

If we have algorithms to compute the local height pairings,
we can solve for the global height pairing in terms of the basis of symmetric bilinear
pairings
on $J(\Q)\otimes\Q_p$ defined by
\begin{equation}\label{E:gij}
g_{ij}(D,E)\colonequals  \frac{1}{2}(\log(D)
(\omega_i)\log(E)(\omega_j) + \log(D)(\omega_j)\log(E)(\omega_i))
\,,\;\, 0\le i\le j\le g-1\,.
\end{equation}
Since we can express $\pi_1(\mathrm{A}(x))$ and $\pi_2(\mathrm{A}(x))$ in terms of the dual basis
$\{\omega_i^*\}$, we can compute
$g_{ij}(\pi(\mathrm{A}(x)))$ for $x\in X(\Q_p)$ (with the obvious abuse of notation)
and extend $h$ to a locally analytic function $h\colon X(\Q_p)\to \Q_p$.

It remains to discuss the computation of the local heights. For  $D_1,D_2\in \Div^0(X_{\Q_p})$ with disjoint
support, the local height is the Coleman integral
certain differential with
residue divisor $\Res(\omega_{D_1})=D_1$ and $c_p$ is a
constant so that 
$c_p^{-1}\chi_p$ extends to a branch $\Q_p^\times \to \Q_p$ of the $p$-adic logarithm; the Coleman integral is taken with
respect to this branch.
The differential $\omega_{D_1}$ is normalised with respect to the splitting
$s$ using

a homomorphism 
$$\Psi\colon T(\Q_p)/T_l(\Q_p)\rightarrow \HH_{\dR}^1(X),$$  from $T(\Q_p)$
the group of differentials of the third kind with integer residues on $X$ quotiented by
$T_l(\Q_p)$ the group of logarithmic differentials $\frac{df}{f}$ with
$f\in \Q_p(X)^*$, as in the algorithm below. We restrict to degree
zero divisors of the form $P-Q$ where $P,Q$ are non-Weierstrass points in
$X(\Q_p)$ that do not reduce to a Weierstrass point in $X(\F_p)$ since we
will need to compute Coleman integrals between $P,Q$, and our 
implementation assumes that these points are in non-Weierstrass disks and
defined over $\Q_p$.

\begin{algorithm}[The local height $h_p(D_1,D_2)$ at $p$ of the global $p$-adic height  {\cite{bb:heights}}]\label{alg:hp}  \hfill \\
Input:
\begin{itemize} \item Hyperelliptic curve $X/\Q_p$, given by an affine
      model $y^2=f(x)$, where $f\in \Z_p[x]$ is squarefree of degree
      $2g+1>2$
\item Prime $p>2g-1$ of good reduction
\item Choice of isotropic subspace $W$ of $H^1_{\dR}(X_{\Q_p})$, complementary to the subspace of regular 1-forms $ \HH^0(X_{\Q_p},\Omega^1)$
\item Divisors $D_1 = P-Q, D_2 = R-S$, where
  $P,Q,R,S$ are non-Weierstrass points in $X(\Q_p)$ that do not reduce to a
    Weierstrass point in $X(\F_p)$, and $R,S$ do not lie in the residue
    disks of $P,Q$.
\end{itemize}
Output: The local height $h_p(D_1,D_2)$ at $p$ of the Coleman--Gross global
  $p$-adic height. 
\begin{enumerate}
	\item Choose $\omega$ a differential in $T(\Q_p)$ with $\text{Res}(\omega)=D_1$.
	
	\item Solve for the coefficients $b_i$ of  $\Psi(\omega)=\sum_{i=0}^{2g-1}b_i\omega_i\in \HH_{\dR}^1(X)$ by
          computing residues, as in \cite[\S 5.2]{bb:heights}. Then $\Psi(\omega)-\sum_{i=0}^{g-1}b_i\omega_i\in W$. Let
	$$\omega_{D_1}\colonequals \omega-\sum_{i=0}^{g-1}b_i\omega_i.$$
      \item
Set $\alpha \colonequals \phi^*(\omega)-p(\omega)$.
        Use Frobenius equivariance of the map $\Psi$ (and the matrix of
        Frobenius computed with respect to the basis $\{\omega_i\}$ of $\HH_{\dR}^1(X)$) to compute
	$$\Psi(\alpha)=\phi^*\Psi(\omega)-p\Psi(\omega)\,.$$
	\item Let $\beta$ be a $1$-form with $\Res(\beta)=(R)-(S)$. Compute $\Psi(\beta)$.
        \item Compute $$h_p(D_1,D_2)\colonequals \int_{D_2}\omega_{D_1} = \int_{S}^R \left(\omega - \sum_{i=0}^{g-1}b_i\omega_i\right),$$
where $$\int_{S}^{R}\omega=\frac{1}{1-p}\left(\Psi(\alpha)\cup\Psi(\beta)+\sum_{A \in X(\C_p)}\textup{Res}_{A}\left(\alpha\int\beta\right)-\int_{\phi(S)}^{S}\omega-\int_{R}^{\phi(R)}\omega\right),$$ see \cite[Remark 4.9]{bb:heights}.
\end{enumerate}
\end{algorithm}
\begin{Remark}Note that in the last step above, $\int_{\phi(S)}^S\omega$
  and $\int_{R}^{\phi(R)}\omega$ are  tiny integrals, that is, Coleman
  integrals between points in the same residue disk. Such integrals may be
  computed merely using a uniformising parameter at any point in the
  residue disk. The computation $\sum_{A \in
  X(\C_p)}\Res_{A}\left(\alpha\int\beta\right)$ will, in most cases,
  require working over various extension of $\Q_p$ to pick up all
  contributions at all poles (see~\cite[Remark~4.10]{bb:heights}).
\end{Remark}
\begin{Remark} If our hyperelliptic curve $X$ does not admit an odd degree model
  over $\Q$, we may choose our prime $p$ such that $X$ has an odd degree model over
  $\Q_p$ and compute local heights at $p$ on this model. This follows from the fact that
  $\Psi(\varphi^*\omega) = \varphi^*(\Psi(\omega))$ 
  for $\varphi$ an isomorphism of curves and
$\omega$  a differential of the third kind.\end{Remark}
\begin{Remark}\label{R:stevan}
  In his thesis~\cite{Gaj22},  Gajovi\'{c} has improved Algorithm~\ref{alg:hp} and extended it to even degree models of hyperelliptic curves.
\end{Remark}

The local height at a prime $\ell \ne p$ is defined in terms of intersection theory. We can
extend $D_1$ and $D_2$ to divisors $\mathcal{D}_1$ and $\mathcal{D}_2$ on a regular model
of $X_{\Q_\ell}$ so that both $\mathcal{D}_i$ have trivial intersection multiplicity with
all vertical divisors;
then by~\cite[Proposition~1.2]{CG89}, we have
$$h_{\ell}(D_1,D_2) = -(\mathcal{D}_1\cdot \mathcal{D}_2)\chi_p(\ell)\,.$$

\subsection{Mordell--Weil sieving}\label{subsec:mws}
The idea of the Mordell--Weil sieve, originally due to Scharaschkin
\cite{Scharaschkin:Thesis}, is to deduce information on rational points on $X$ via the
intersection of the images of $X(\F_v)$ and $J(\Q)$ in $J(\F_v)$ (or suitable quotients) for several primes $v$ of
good reduction.
It is often applied to verify that $X(\Q)=\varnothing$, but it can also be
combined with
$p$-adic techniques to compute $X(\Q)$ when there are rational points.

We review the basic idea, which is straightforward. Making the sieve perform well in practice
is a different  matter; see~\cite{Bruin-Stoll:MWSieve} for an elaborate discussion of the
issues one encounters and detailed strategies.
%solutions. 
For ease of exposition, we assume that $J(\Q)$ is torsion-free and that we have
generators $P_1,\ldots,P_r$ of $J(\Q)$. 
Let $M>1$ be an integer and let $S$ be a finite set of primes of good reduction for $X$.
Then the diagram
\[
    \xymatrix{
    X(\Q)\ar@{->}[r]  \ar@{->}[d]&
    J(\Q)/MJ(\Q)\ar@{->}[d]^{\alpha_{S,M}}\\
    \prod_{v\in S}X(\F_v)\ar@{->}[r]_-{\beta_{S,M}}&\prod_{v\in S} J(\F_v)/MJ(\F_v)\,\\
} \]
is commutative. In the situation of interest to us, the horizontal maps are induced by our
choice of base point $b \in X(\Q)$. 

In our work, we use the Mordell--Weil sieve in two ways.
On the one hand, we apply it to show that for a fixed prime $p$, a given
residue disk in $X(\Q_p)$ does not contain a rational point. To this end,  we 
set $M = M'\cdot p$ for some suitable auxiliary integer $M'$, 
and we choose $S$ to consist of primes $\ell$ so that $\gcd(\#J(\F_{\ell}), \#J(\F_q))$ is
large for some prime divisors $q\mid pM'$.
We can then hope that the image of the reduction of the disk under $\prod
\beta_{S,M}$ does not meet the
image of the map $\prod \alpha_{S,M}$.

On the other hand, we use the sieve to show for fixed $M>1$ that a given coset of $MJ(\Q)$ does not contain
the image of a point in $X(\Q)$ under the Abel--Jacobi map $P\mapsto [P-b]$.
Suppose a point $P\in X(\Q_p)$ is given to finite precision $p^N$.
If $P$ is rational, then there are integers $a_1,\ldots,a_g$ such that 
\[
  [P-b]= a_1P_1+\cdots+a_gP_g.
\]
Via the abelian logarithm, we compute a tuple
$(\tilde{a}_1,\ldots,\tilde{a}_g) \in \Z/p^N\Z$  satisfying $a_i \equiv \tilde{a}_i \pmod{p^N}$ for all
$i\in \{1,\ldots,g\}$. To show that $P$ is not rational, it suffices to show 
that the corresponding coset of $p^NJ(\Q)$ does not contain the image of such a point.

In our implementation, we have not tried to optimise the interplay between quadratic
Chabauty and the Mordell--Weil sieve. Such an optimisation is discussed
in~\cite[\S7]{BBM17}.
Let us only note here that we may combine quadratic Chabauty information coming from several
primes, and that we can enhance that information using an auxiliary integer $M'$ similar
to the above.
Another account of combining quadratic Chabauty with the Mordell--Weil sieve
can be found in~\cite[\S 6.7]{BBBLMTV19}.

\begin{Remark}
  All examples in this paper 
  
  satisfy $r=g=\rk_{\Z}\mathrm{NS}(J)$, resulting in at least two independent locally
analytic functions vanishing in $X(\Q)$ for the $g>2$ examples. Since we expect that their common zero set is precisely $X(\Q)$ (or that there is a geometric reason for the appearance of any additional $p$-adic solutions), we do not expect to require the sieve.
Indeed, we only had to apply the sieve for curves of genus~2.
For these examples, we always required only one prime for the quadratic Chabauty
computation; we chose this prime in such a way as to simplify the sieving.
\end{Remark}

\subsection{Implementation and scope}\label{subsec:imp}

We have implemented the algorithms described in this section in the computer algebra system
{\tt Magma}~\cite{BCP97}. Our code is freely available at~\cite{QCCode}. It
extends the code used for $X_{\mathrm{s}}^+(13)$ in~\cite{BDMTV19} and can
be used to recover that example. It is applied to new examples, as
discussed in \S\ref{sec:examples}. 

\par We begin by summarising our discussion so far and describe the general procedure to
determine a finite set containing $X(\Q_p)_2$ as it would apply to the modular curve $X$ attached
to a general congruence subgroup, and Atkin--Lehner quotients thereof. In this generality,
several steps cannot be easily automated, so we discuss the extent to which our
implementation has automated the procedure, and point out which steps require additional
action from the user. See Example~\ref{E:107} for a fairly detailed worked example.

Our techniques are built on prior work of Tuitman on computing the action of Frobenius on rigid cohomology \cite{Tui17}. We recall some of the underlying structures present in Tuitman's work  and a set of assumptions on these auxiliary structures. 

Suppose our modular curve $X/\Q$ is given by a (possibly singular)  plane
model $Q=0$ with $Q(x,y) \in \ZZ[x,y]$ a polynomial that is irreducible and
monic in $y$. Let $d_x$ and $d_y$ denote the degrees of the morphisms $x$
and $y$, respectively, from $X$ to the projective line. Let $\Delta(x) \in
\ZZ[x]$ denote the discriminant of $Q$ with respect to the variable $y$.
Moreover, define $r(x) \in \ZZ[x]$ to be the squarefree polynomial with the
same zeroes as $\Delta(x)$,
in other words, $r=\Delta/(\gcd(\Delta,\frac{d\Delta}{dx}))$.

\begin{Definition}\label{defn:intbases} Let $W^0 \in \GL_{d_x}(\QQ[x,1/r])$ and
$W^{\infty} \in \GL_{d_x}(\QQ[x,1/x,1/r])$ denote matrices such that, if we denote 
\[ b^0_j = \sum_{i=0}^{d_x-1} W^0_{i+1,j+1} y^i \; \; \; \; \mbox{ and } \; \; \; \; b^{\infty}_j = \sum_{i=0}^{d_x-1} W^{\infty}_{i+1,j+1} y^i \] 
for all $0 \leq j \leq d_x-1$, then
\begin{enumerate}
\item $[b^{0 \;}_0,\ldots,b^{0 \;}_{d_x-1}]$ is an integral basis for $\QQ(X)$ over $\QQ[x]$,
\item $[b^{\infty}_0,\ldots,b^{\infty}_{d_x-1}]$ is an integral basis for $\QQ(X)$ over $\QQ[1/x]$,
\end{enumerate}
where $\QQ(X)$ denotes the function field of $X$. Moreover, let $W \in \GL_{d_x}(\QQ[x,1/x])$ denote
the change of basis matrix $W=(W^0)^{-1} W^{\infty}$.
  
\end{Definition}

\begin{assump}[{\cite[Assumption 1]{Tui17}}]\label{tuitman1} \mbox{ }
\begin{enumerate}
\item The discriminant of $r(x)$ is contained in $\ZZ_p^{\times}$.
\item If we denote 
$b^0_j = \sum_{i=0}^{d_x-1} W^0_{i+1,j+1} y^i$ and 
$b^{\infty}_j = \sum_{i=0}^{d_x-1} W^{\infty}_{i+1,j+1} y^i$ for all $0 \leq j \leq d_x-1$, and if we 
let $\FF_p(x,y)$ be the field of fractions of $\FF_p[x,y]/(Q)$, then:
\begin{enumerate}
\item The reduction modulo~$p$ of $[b^{0 \;}_0,\ldots,b^{0 \;}_{d_x-1}]$ is an integral basis for $\FF_p(x,y)$ over $\FF_p[x]$.
\item The reduction modulo~$p$ of $[b^{\infty}_0,\ldots,b^{\infty}_{d_x-1}]$ is an integral basis for $\FF_p(x,y)$ over $\FF_p[1/x]$.
\end{enumerate}
\item $W^0 \in \GL_{d_x}(\ZZ_p[x,1/r])$ and $W^{\infty} \in \GL_{d_x}(\ZZ_p[x,1/x,1/r])$.
\item Denote:
\begin{align*}
\mathcal{R}^0        &= \ZZ_p[x]b^{0}_0 \; \; \; \; \; +\ldots+\ZZ_p[x]b^{0}_{d_x-1}, \\
\mathcal{R}^{\infty} &= \ZZ_p[1/x]b^{\infty}_0+\ldots+\ZZ_p[1/x]b^{\infty}_{d_x-1}.
\end{align*}
For a ring $R$, let $R_{\textrm{red}}$ denote the reduced ring obtained by quotienting out by the nilradical. Then the discriminants of the finite $\ZZ_p$-algebras $(\mathcal{R}^0/(r(x)))_{\textrm{red}}$ and
$(\mathcal{R}^{\infty}/(1/x))_{\textrm{red}}$ are contained in $\ZZ_p^{\times}$.\end{enumerate}
\end{assump}

\begin{Remark} 
These conditions imply that the curve $X$ has good reduction at~$p$.
\end{Remark}

\begin{algorithm}[Quadratic Chabauty for modular curves]\label{alg:qc} $\qquad$ \;\\
Input:  
\begin{itemize}
\item A modular curve $X/\Q$ with Mordell--Weil rank $r=g$ and $\rk_{\Z}\mathrm{NS}(J)>1$, and for which the image of $J(\Q )$ in $\HH^0
  (X_{\Q _p },\Omega^1 )^\vee$ has rank $g$.
  \item  A covering of $X$ by affine opens that are birational to a planar curve cut out by an equation that is monic in one variable, has $p$-integral coefficients and satisfies Assumption \ref{tuitman1}. (See \S\ref{subsubsec:comp-patches}.)

\item A prime $p$ of good reduction such that the Hecke operator $T_p$ generates $\End ^0 (J)$.
\item   For all primes $\ell $ that are not of potentially good reduction, the local height functions $X(\Q _{\ell })\to \Ker (\NS (J)\to \NS (X))_{\Q _p }^*$, computed using Theorem \ref{thm:BeD}. (See  \S\ref{subsubsec:comp-heights-away}.)

    \item A starting precision $n$.   
        \item A height bound $B$.
\end{itemize}

Output:  An approximation to a finite set containing the set of points $X(\Q_p)_2$, computed to precision $n' \leq n$ or FAIL. 

\begin{enumerate}
  \item\label{smallpts} Compute the set $X(\Q)_{\mathrm{known}}$ of points  in $X(\Q)$
    with height bounded by $B$.
  \item\label{intsympbasis} Compute an integral symplectic basis for
    $\HH^1_{\dR}(X_{\Q})$ or return FAIL. 
  
  \item\label{frobhecke} Compute the action of Frobenius on
    $\HH^1_{\dR}(X_{\Q_{p}})$ using Tuitman's algorithm~\cite{Tui16,Tui17}. Use the Eichler--Shimura relation to compute the matrix of the action of the Hecke operator $T_{p}$ on $\HH^1_{\dR}(X_{\Q_{{p}}})$.
  \item\label{hodgeequiv} Compute a splitting of the Hodge filtration that is equivariant for the action of $\mathrm{End}(J)$ in the sense of \S\ref{subsec:global_height}.   \item\label{corrs} Compute the matrices of a basis $Z_1, \ldots, Z_{\rk\NS(J) - 1}$ of
  $\mathrm{Ker}\left(\mathrm{NS}(J) \to \mathrm{NS}(X)\right)$ acting on $\HH^1_{\dR}(X_{\Q_{p}})$, see \S\ref{subsubsec:comp-NS}. 
\item\label{inita} Let $A\colonequals \varnothing$. For each $Z_i$, compute the associated heights: 
\begin{enumerate}

  \item\label{aff} For each affine patch, do the following:
\begin{enumerate}

  \item\label{lambdafil} Compute the functions $\lambda ^{\Fil }$ from \eqref{eqn:definitions-abcs} using~\cite[\S 4]{BDMTV19}.
\item\label{lambdaphi} Compute the functions $\lambda ^{\phi }$ from \eqref{eqn:definitions-abcs} using~\cite[\S 5]{BDMTV19}.
\end{enumerate}
\item\label{solve} Solve for the height pairing, either using a large enough supply of
  known rational points $P_1, \ldots, P_n$ on $X$, if possible, or by
    computing the Coleman--Gross height pairing on $r$ independent points
    in $J(\Q)$. (See \S\ref{subsubsec:comp-globalheight}.) If this is
    unsuccessful, return FAIL. 

\item\label{roots} Compute solutions of the function(s) coming from $Z_i$ or return FAIL
  if there has been too much precision loss to determine these solutions. 
\item\label{simple} Check that the solutions are simple. If there is a non-simple
  solution corresponding to a point in  $X(\Q)_{\mathrm{known}}$, return FAIL. Else, add to
    the set $A$ the solutions that (simultaneously) satisfy the(se) function(s).
\end{enumerate}
\item\label{return} Return $A$. 
\end{enumerate}
\end{algorithm}
\begin{Remark}
We assume that we know a priori that the Mordell-Weil rank of the Jacobian is equal to the genus
  of the curve. For modular curves, by Gross--Zagier--Kolyvagin--Logachev
  this amounts to checking that the associated eigenforms have analytic
  rank one (see e.g.~\cite[\S 7]{DLF19}). For hyperelliptic curves, it is
  sometimes simpler to carry out a two-descent. 
\end{Remark}
\begin{Remark}
Note that if the algorithm fails due to a loss of precision, it may be
  possible to remedy this by increasing the starting precision. One place
  where increasing precision may not work is if the $p$-adic logarithm does
  not induce an isomorphism $J(\Q )\otimes \Q _p \simeq \HH^0 (X_{\Q _p
  },\Omega^1 )^\vee $, even though the rank of $J(\Q)$ is $g$. For the Atkin-Lehner quotients $X_0 ^+ (N)$, the weak Birch--Swinnerton-Dyer conjecture implies $J(\Q )$ always generates $\HH^0 (X_{\Q _p },\Omega^1 )^\vee $ (see \cite[Lemma 7]{DLF19}). In general, if $r=g$ and the Zariski closure of $J(\Q )$ is $J$, then a conjecture of Waldschmidt \cite[Conjecture 1]{Waldschmidt} (an analogue of the Leopoldt conjecture for abelian varieties) implies that the $p$-adic logarithm is always an isomorphism. In theory, if one \textit{knew} that $J$ gave a counterexample to Waldschmidt's conjecture, and $r=g$, then one could simply apply the Chabauty--Coleman method. However, a priori it could happen that $J$ gave a counterexample but there was no way of verifying this by a computation to finite $p$-adic precision.
  Another place where increasing precision will not help is if there are multiple
  roots in Step~\eqref{roots}. However, we only expect this to happen for geometric
  reasons. 
  
One can have $r>g$ for the curves $X_0 ^+ (N)$ with $N$ prime, even though $X_0 ^+ (N)(\Q _p )_2 $ is always finite when the genus is greater than one \cite{DLF19}. However the smallest genus for which this happens is $g=206$ (with $N=5077$), so the $r=g$ hypothesis is not the main restriction to the scope of our algorithms for this family of curves.
\end{Remark}

\begin{Remark}\label{rem:unitroot}In the case when $p$ is a prime of ordinary reduction for the Jacobian, one may take the splitting of the Hodge filtration given by the \emph{unit root subspace}, that is, the unit root eigenspace of Frobenius $\phi$ acting on $\HH^1(X_{\Q_p})$. Given a basis $\{\eta_1, \ldots, \eta_{2g}\}$ of $\HH^1(X_{\Q_p})$, where $\eta_1, \ldots, \eta_g$ are holomorphic, a basis for the unit root eigenspace mod $p^n$ is given by $\{(\phi^*)^n\eta_{g+1},\ldots, (\phi^*)^n\eta_{2g}\}.$ \end{Remark}
\begin{Remark}
In this paper, we do not discuss algorithms for computing the input of the local height
  functions as maps from $\Q _\ell$-points to $\Q_p$-linear functionals on $\Ker (\NS (J) \to \NS (X))$. In
  Section \ref{sec:examples} we give examples where this function can be nontrivial, and
  where $X(\Q )$ can still be determined using quadratic Chabauty. There are two
  procedures we illustrate for doing this. In Section \ref{subsec:nontriv_away}, we
  calculate regular semistable models at bad primes and have a sufficient supply of
  rational points (and sufficiently simple dual graphs) to reconstruct the
  functions $j_\ell$
  from Theorem \ref{thm:BeD} using evaluation of $p$-adic local heights at known rational points. In Section \ref{sec:ns17}, although we know a regular semistable model ``abstractly,'' we do not know the relation between the stable model (at the bad prime 17) and the model we use for $p$-adic calculations. This, together with the relative paucity of known rational points, makes it infeasible to apply the first procedure. Instead, we use extra information about the action of inertia on the stable model, together with Theorem \ref{thm:BeD} to identify a subspace of line bundles in $\Ker (\NS (J)\to \NS (X))$ for which the associated local heights vanish.
\end{Remark}
To further determine the subset of rational points $X(\Q)$ from the finite set of
points produced by our algorithm, we carry out the Mordell--Weil sieve. In practice it may
happen (see below) that $X(\Q)$ is returned by the algorithm, but this is typically not
the case when $X$ has genus two.

\subsubsection{Affine patches} 
\label{subsubsec:comp-patches}
Most of the examples discussed in Section \ref{sec:examples} are either hyperelliptic curves or smooth plane quartics.  As demonstrated in Section \ref{sec:ns17}, our code is sometimes able to treat more general examples. Our implementation was designed to take as input a plane affine patch $Y\colon Q(x,y)=0$ of a modular curve $X/\Q$ satisfying the requirements in~\S\ref{S:setup} and a prime $p$ of good reduction. It returns all rational points on $X$ in affine residue disks where the lift of Frobenius constructed in~\cite{Tui16,Tui17} is defined. Note that we do not require $Y$ to be smooth, but we need $Q$ to be monic with $p$-integral coefficients. 

\par We can sometimes find an affine patch $Y$ having the convenient property that all rational points on $X$ must be among the points returned by running our algorithm on $Y$. If no such $Y$ is found, then we need to find two suitable affine patches such that every rational point on $X$ is contained in at least one patch. For smooth plane quartics, our implementation includes an algorithm that automates this process for convenience of the user. For other curves, this step is left to the user.

\subsubsection{The N\'eron--Severi classes $Z_i$}
\label{subsubsec:comp-NS}

Under the assumption that $T_p$ generates the endomorphism ring of the Jacobian, which we
made for convenience above, one may proceed precisely as in~\cite[\S6.4]{BDMTV19} to
determine a nontrivial class 
\[
Z \ \in \ \mathrm{Ker}\left(\mathrm{NS}(J) \longrightarrow \mathrm{NS}(X) \right).
\]
Indeed, the matrix $A_p$ of the Hecke operator $T_p$ acting on $\HH^1_{\dR}(X_{\Q_p})$ is easily determined from the matrix of Frobenius $F_p$ (which is already a byproduct of the algorithms for the local height at $p$), by the Eichler--Shimura relation: $$T_p = F_p +pF_p^{-1}.$$

Under our assumption, the matrices of the classes $Z_i$ acting on $\HH^1_{\dR}(X_{\Q_p})$ may then be computed as linear combinations of powers of $A_p$.

\begin{Remark}

This is the only part of our algorithm specific to modular curves, since it relies on the
  Eichler--Shimura relation. It should however be noted that this is mainly a matter of
  convenience adopted for the purpose of automation. More generally, for a smooth
  projective curve $X/\Q$ satisfying the assumptions of \S\ref{S:setup}, one could find
  $p$-adic approximations of the action of the nontrivial classes $Z_i$ on
  $\HH^1_{\dR}(X_{\Q_p})$ using just $p$-adic linear algebra. Indeed, the
  space of correspondences which are symmetric under the Rosati involution and induce endomorphisms of trace zero on the Tate module 
  maps under the cycle class into the intersection of the $\Fil ^{1}$ and $\phi =p$ subspaces of 
\begin{equation}\label{eqn:BO}
\ker \left( \wedge ^2 \HH ^1 _{\dR}(X_{\Q _p } )\stackrel{\cup }{\lra } \HH ^2 _{\dR}(X_{\Q _p } )
     \right)\, .
\end{equation}
In fact, by the Tate conjecture, the rank of the space of (crystalline)
  cohomology classes of such correspondences over $\F _p $ is equal to the
  dimension of the $\phi =p$ subspace of \eqref{eqn:BO}, and by the
  $p$-adic Lefschetz-(1,1) theorem of Berthelot and Ogus \cite[\S
  3.8]{BO83} 
  such a correspondence over $\F _p$ lifts to $\Q _p$ if and only if 
  its cycle class lies in $\Fil ^{1}$. Note that the dimension of the space of correspondences symmetric under the Rosati involution need not equal the dimension of $\wedge ^2 \HH ^1 _{\dR}(X_{\Q _p })^{\phi =p}\cap \Fil ^1$, as was erroneously claimed in \cite[Lemma 4.5]{BDMTV19}, since the rank of the intersection of a $\mathbf{Z}$-lattice with a $\Q _p$-subspace may be less than the dimension of the intersection with the $\Q _p$-subspace it spans. However, if one knows a set of generators of a finite index subgroup of $\End (J)$ in advance (e.g. using algorithms for rigorous computation of the endomorphism algebra of the Jacobian \cite{CMSV})) then one can use this to compute the classes of generators in cohomology. 

Therefore the assumption that $T_p$ generates the endomorphism algebra could be
  circumvented in this step with a little work, although it is used in the computation of
  the local heights away from $p$, see below. When the assumption is not satisfied, our
  implementation throws an error, urging the user to try a different choice of prime $p$.

\end{Remark}

\subsubsection{The local heights away from $p$}
\label{subsubsec:comp-heights-away}

This step requires an explicit knowledge of a semi-stable model of the modular curve $X$, as well as a description of the action of $Z_i$ on the concomitant cohomological structures in order to be able to apply Theorem \ref{thm:BeD}. It is clear that a full automation of this step, starting from a set of defining equations for $X$, falls outside the scope of our implementation.

\par Semi-stable models for modular curves are known in many cases, see for
instance the recent work of Edixhoven--Parent \cite{EP21}. In practice, one
can also often use the \texttt{SageMath} toolbox {\tt MCLF} \footnote{
  {\url{https://github.com/MCLF/mclf}}} due to R\"uth and Wewers to compute
such models. The main advantage of having computed the $Z_i$ in
\S\ref{subsubsec:comp-NS} as combinations of powers of $T_p$ is that this
makes it easier to compute the quantities appearing in Theorem
\ref{thm:BeD}. Even though we see no way to fully automate this step, we
hope to convince the reader of its practicality by working it out for
the genus~2 curves $C_{188}$ and $C_{161}$ in Examples~\ref{E:188}
and~\ref{E:161}. 

\subsubsection{The global height pairing}
\label{subsubsec:comp-globalheight}
If there are not sufficiently many rational points on the curve to solve for the height
pairing, we instead compute the local heights $h_v$ in the sense of Coleman and
Gross, see \S\ref{subsec:coleman-gross}. 
For hyperelliptic curves $X/\Q_p$ of odd degree, $h_p(D_1,D_2)$ can be computed using an algorithm
due to Balakrishnan--Besser  \cite{bb:heights, bb:heightserrata}. Based on earlier
{\tt SageMath} code due to Balakrishnan, we have implemented this in {\tt
Magma} for divisors $D_1$ and $D_2$ that split over
$\Q_p$, have support contained in disjoint residue disks and for which no points in the support reduce to Weierstrass points mod $p$. 
To compute the local heights $h_\ell$ for $\ell\ne p$, we rely on {\tt Magma}'s implementation of
an algorithm for local canonical heights on hyperelliptic curves described by Holmes and M\"uller~\cite{Hol14,
Mul14}. An algorithm for general curves was given by van Bommel, Holmes and M\"uller \cite{van2020explicit}.

To solve for the height pairing, we need
to find representatives for $r$ independent points in $J(\Q)$  that satisfy the
 assumptions mentioned above.
Our implementation is currently restricted to genus 2 curves, since this
step was only necessary for such curves, 
but a generalisation to higher
genus hyperelliptic curves would be straightforward.

\begin{Remark}
The code is currently restricted to the base field $K=\Q$. To extend it to
  more general number fields, one would need to combine these algorithms
  with those used in \cite{BD18} for imaginary quadratic fields in certain
  cases, or for general number fields, with those in \cite{BBBM}.
\end{Remark}

%%=========================================================================
\section{Precision analysis}\label{sec:prec}
In this section, we bound the loss of absolute $p$-adic precision that may occur in
our computations by bounding the valuations of the error terms. We 
also estimate the valuations of the power series expansion of the quadratic Chabauty
function $\rho$ and use this to bound the precision of its roots.

We keep the notation used in the previous sections.
Recall from~\eqref{eq:rho} that $\rho = h-h_p$, where
\begin{itemize}
  \item $h$ is the global $p$-adic height defined in~\eqref{eq:hdef};
  \item $h_p$ is the local component of $h$, discussed
    in~\S\ref{subsec:local_decomp}.
\end{itemize}
By~\eqref{hpformula}, the local height $h_p$ satisfies
\[
  h_p(x) = \gamma_{\phi}
 - \gamma_{\Fil}
 -\boldsymbol{\beta}^\intercal_{\phi} \cdot s_1 (\boldsymbol{\alpha}_{\phi})
 - \boldsymbol{\beta }^\intercal_{\Fil} \cdot s_2
  (\boldsymbol{\alpha}_{\phi})\,,
\]
where the Hodge filtration of the filtered
$\phi$-module
$M(x) \colonequals \left( \mathrm{A}_Z(b,x) \otimes_{\Q_p}
\mathrm{B}_{\cris}\right)^{G_{\Q_p}} $ discussed
in~\S\ref{subsec:local_decomp} is encoded by
$\boldsymbol{\beta}_{\Fil}$ and $
\gamma_{\Fil}$ 
and $\boldsymbol{\alpha}_{\phi}, \boldsymbol{\beta}_{\phi}$ and $
{\gamma}_{\phi}$ encode the Frobenius structure of $M(x)$.

We will bound the loss of precision in the computation of the Hodge
filtration in~\S\ref{subsec:hodge}, and we do the same for the Frobenius
structure in~\S\ref{subsec:frobsplit}.
In~\S\ref{subsec:globht}, we bound the precision loss for the global height
computation. In the final part of this section,~\S\ref{subsec:boundcoeffs} we
bound the valuation of the coefficients of the expansion of $\rho$ in a
residue disk, and we discuss how this may be used to provably determine the
roots of $\rho$ to a certain precision.
This section relies heavily on~\cite[Sections 4,5]{BDMTV19}.

\subsection{Hodge filtration}\label{subsec:hodge} 
We first bound the loss of precision in
Steps~\eqref{intsympbasis}--\eqref{corrs} of Algorithm~\ref{alg:qc}. For
simplicity, we restrict to one class $Z$; the extension to $\rk\NS(J)
- 1$ such classes is immediate.
Let $Y/\Q$ be an affine open subset of $X$, birational to a curve given by an equation
that satisfies Assumption~\ref{tuitman1}. We may compute an integral, symplectic basis $\boldsymbol{\omega }=(\omega _0 ,\ldots ,\omega _{2g-1})$ of de Rham cohomology over $\Q $ exactly, and extend this to an integral basis of $H^1 _{\dR}(Y)$ via differentials $(\omega _{2g},\ldots ,\omega _{2g+d-2})$ of the third kind. 
Using such a basis, we may compute the action of the Frobenius operator $F$
on $\HH ^1 _{\dR}(X/\Q _p )$ to any desired $p$-adic precision using
Tuitman's algorithm~\cite{Tui16, Tui17}, from which we obtain the action of the Hecke operator $T_p = F + pF^{-1}$ on $\HH^1 _{\dR}(X/\Q _p)$ by Eichler--Shimura. The inversion of $F$ results in a finite and computable loss of precision, which the code takes into account. This results in an algorithm that returns the action of the correspondence $Z$ correctly modulo $p^n$ for some $n\geq 1$ that is returned by the algorithm.  

Using this, we may compute a matrix $\Lambda $ with entries in $H^0 (Y,\Omega _{Y_{\Q}})$, of the form
$$\Lambda \colonequals  -\left( \begin{array}{ccc}
0 & 0 & 0 \\
\boldsymbol{\omega} & 0 & 0 \\
\eta & \boldsymbol{\omega}^{\intercal}Z & 0 \\
\end{array} \right)\,$$
such that
$d+\Lambda$ extends to a flat connection on $X$. From this, we may compute $\gamma _{\Fil }$ and $\boldsymbol{\beta}_{\Fil }$ from \eqref{eqn:definitions-abcs}. We recall from \cite[\S 4]{BDMTV19} that the defining properties of $\eta $, the $\boldsymbol{\beta}_{\Fil }$ and $\gamma _{\Fil}$ are as enumerated below. For $x\in (X-Y)(\overline{\Q })$, we let $t_x$ denote a parameter, and $\boldsymbol{\Omega }_x $ denote the vector of formal integrals of the basis differentials $\omega _i $:
\[
d\boldsymbol{\Omega }_{x,i}=\omega _i \in\overline{\Q }[\! [t_x]\! ].
\]
\begin{enumerate}
\item The first $g$ entries of $\boldsymbol{\beta }_{\Fil }$ are zero, and the last $g$ are given by a vector $\mathbf{b}_{\Fil }$ of constants specified below.
\item $\eta $ is  
  a linear combination of $\omega _{2g},\ldots ,\omega _{2g+d-2}$,
    unique by~\cite[Lemma~4.10]{BDMTV19} such that
\begin{equation}\label{eqn:hodge_defn1}
d\boldsymbol{\Omega }_x^{\intercal }Z \boldsymbol{\Omega }_x -\eta
\end{equation}
has vanishing residues at all $x\in (X-Y)(\overline{\Q })$.
\item 
$\mathbf{b}_{\Fil }$ and $\gamma _{\Fil } \in \mathcal{O}(Y)$ are the unique solutions to the equation $\gamma _{\Fil }(b)=0$ and 
\begin{equation}\label{eqn:hodge_defn2}
g_x +\gamma _{\Fil}-\mathbf{b}_{\Fil }^\intercal N^\intercal \boldsymbol{\Omega }_x -\boldsymbol{\Omega }_x ^\intercal Z N N^\intercal \boldsymbol{\Omega }_x \in L[\! [t_x ]\! ]
\end{equation}
for all $x\in (X-Y)(\overline{\Q })$, where $g_x \in \overline{\Q }[\! [t_x ]\! ]$ is defined to be the formal integral of $d\boldsymbol{\Omega }_x ^\intercal Zd\boldsymbol{\Omega }_x -\eta $ and $N$ is the block $2g \times g$ matrix with top block zero and lower block a $g \times g$ identity matrix.
 
\end{enumerate}
Given our basis $\boldsymbol{\omega }$, we may calculate $\boldsymbol{\Omega }_x $ to any given $t_x$-adic precision. Note that to solve \eqref{eqn:hodge_defn1}, we only need to know $\boldsymbol{\Omega }_x$ modulo $t_x ^{m_x}$, where $m_x$ is the maximum of the order of the poles of the entries of $\boldsymbol{\Omega }_x$. Similarly, to solve for $\gamma _{\Fil}$ and $\mathbf{b}_{\Fil}$ in \eqref{eqn:hodge_defn2}, we need only compute the principal parts of $\boldsymbol{\Omega }_x$ and $\boldsymbol{\Omega }_x ^\intercal Z N N^\intercal \boldsymbol{\Omega }_x $. Hence given the above we may calculate $\eta ,\gamma _{\Fil }$  and $\mathbf{b}_{\Fil }$ to precision $p^{n-2\nu }$, where $\nu $ is minus the minimum of the valuations of the $t_x^i $ coefficients of the entries of $\boldsymbol{\Omega }_x $, for $i\leq m_x$.

\subsection{Frobenius-equivariant splitting}\label{subsec:frobsplit}
We now bound the loss of precision in the computation of the Frobenius-equivariant splitting
$$
\lambda^\phi(x) =\left(
\begin{array}{ccc} 
1 & 0 & 0 \\
  \boldsymbol{\alpha }_{\phi}(b,x) & 1 & 0 \\
  \gamma_{\phi}(b,x) & \boldsymbol{\beta }^{\intercal}_{\phi}(b,x) & 1 \\
\end{array}
\right)$$
from~\eqref{eqn:definitions-abcs} for $x \in X(\Q_p)\cap ]\mathcal{U}[$
where $\mathcal{U}$ is an open of $Y_{\F _p }$ on which we
have an overconvergent lift of Frobenius. This computation is the content
of~\cite[\S5]{BDMTV19}.

The first step is to find the Frobenius structure on the filtered
$\phi$-module $M(b)$. By~\cite[\S5.3.2]{BDMTV19}, the inverse of the
Frobenius structure is given by a
matrix $$G \in (\HH^0(]Y[,j^\dagger
\mathcal{O}_Y))^{(2g+2)\times(2g+2)}$$ such that
\begin{equation}\label{lambdaG}
  \Lambda_{\phi}G + dG = G\Lambda,
\end{equation}
where $j^\dagger \mathcal{O}_Y$ is the overconvergent structure sheaf on the tube $]Y[$.

Compared to~\cite[\S5.3.2]{BDMTV19}, we give a slightly more detailed
account of the algorithm to find $G$. We first apply the algorithms in~\cite{Tui16,Tui17} (see~\cite[Algorithm~2.18]{BT19}) to
compute the action of Frobenius on
$\HH^1_{\rig}(X\otimes \Q_p)$ as
\begin{equation}\label{}
\phi^*\boldsymbol{\omega} = F\boldsymbol{\omega} + d\mathbf{f} 
\end{equation}
for a matrix $F \in M_{2g}(\Q_p)$ and a column vector $\mathbf{f}$ with entries in
$\HH^0(]Y[,j^\dagger \mathcal{O}_Y)$, uniquely determined by the condition that
$\mathbf{f}(b_0 ) = \mathbf{0}$, 
where $b_0$ is the Teichm\"uller point in the disk of $b$.

Next, we define a vector of functions 
$ \mathbf{g}_0   \colonequals -F^TZ  \mathbf{f}  $.
Then, the differential
\begin{equation}\label{chiexpansion}
  \xi \colonequals(\phi^*\boldsymbol{\omega}^{T})Z\mathbf{f} + (\phi^*\eta - p \eta) -\mathbf{g}_0 ^{T} \boldsymbol{\omega}
\end{equation}

is of the second kind, and therefore the reduction algorithms in
$\HH^1_{\mathrm{rig}}(Y)$ from~\cite{Tui16,Tui17} can be applied to compute a vector of
constants $\mathbf{c} \in \Q_p^{2g}$ and a function $H\footnote{The
function $\h$ is denoted $h$ in~\cite{BDMTV19}, but we chose a different
notation to avoid confusion with the global height, which is also denoted
by $h$.} \in \HH^0(]Y[,j^\dagger
\mathcal{O}_Y)$ such that
  \begin{equation}\label{}
    \mathbf{c}^T \boldsymbol{\omega} + d\h  = \xi.
  \end{equation}

Hence the function $\mathbf{g}\colonequals \mathbf{g}_0 +\mathbf{c}$ satisfies
$$d\mathbf{g}^\intercal  =  d\mathbf{f}^\intercal ZF\quad\text{and }\quad
d\h  =  \boldsymbol{\omega}^\intercal F^\intercal Z\mathbf{f} + d\mathbf{f}^\intercal Z\mathbf{f}-\mathbf{g}^\intercal \boldsymbol{\omega}
+\phi^*\eta - p\eta\,,$$
and we normalise $\h$ by requiring that $\h(b_0)=0$.
The matrix
\begin{equation}\label{eqn:F-structure}
G = \left( \begin{array}{ccc}
1 & 0 & 0 \\
\mathbf{f} & F & 0\\
\h & \mathbf{g}^\intercal & p \\
\end{array} \right)\end{equation}
then satisfies~\eqref{lambdaG}.

\subsubsection{Frobenius-equivariant splitting for Teichm\"uller points}\label{subsec:teich}

Suppose that $x_0\in X(\Q_p)\cap]\mathcal{U}[$ is a Teichm\"uller point.
As described in~\cite[\S5.3.2]{BDMTV19}, the Frobenius-equivariant splitting of $M(x_0)$ is given by
\begin{equation}\label{}
\lambda^\phi(x_0) =
\left( \begin{array}{ccc}
1 & 0 & 0 \\
(I-F)^{-1}\mathbf{f} & 1 & 0\\
\frac{1}{1-p}\left(\mathbf{g}^{\intercal}(I-F)^{-1}\mathbf{f} +\h\right) & \mathbf{g}^{\intercal}(F-p)^{-1} & 1 \\
\end{array} \right)(x_0)\,.
\end{equation}
  The loss of precision in the computation of $\mathbf{f}$ and $F$ is estimated in~\cite{Tui17}. 
Hence it is easy to bound the precision loss in the computation of
  $\lambda^\phi(x_0)$ using the
  following result.
\begin{prop}\label{P:G}
Suppose that the entries of the matrix $G$ and a point $P\in
  X(\Q_p)\cap ]\mathcal{U}[$ are accurate to $n$ digits of precision.  Then $G(P)$ is also accurate
  to $n$ digits of precision. 
\end{prop}

  Our proof of Proposition~\ref{P:G} is somewhat
  similar, but more involved than the proofs in~\cite[\S4]{BT19}, where the
  loss of precision in the evaluation of~$\mathbf{f}$ and of single Coleman
  integrals is estimated.
  We may expand \begin{equation}\label{chi-sum}
    \xi = \sum_{j \in \Z}\left(\sum_{k=0}^{d_x - 1} \frac{w_{j,k}(x)}{r(x)^j}
  b_k^0\right)\frac{dx}{r}\,. 
  \end{equation}
  The hardest part of the proof of Proposition~\ref{P:G} is to
  find lower bounds on the
  valuation of the coefficients $w_{j,k}$, which we now describe.   Let $e_0$ (resp., $e_\infty$) be the maximum of the ramification indices of the map $x\colon X\to
  \mathbf{P}^1$ with respect to our chosen model at points lying in affine (resp., infinite) disks.
\begin{lemma}\label{L:chi-bds}
  There is a constant  $\kappa$ such that for all $j,k$ we have
\begin{equation}\label{ord-omjk}
  \ord_p(w_{jk}) \geq \begin{cases}
&\left\lfloor\frac{j}{p}\right\rfloor + 1 - \log_p(je_0) +\kappa , \quad j \neq 0\\
&\kappa, \qquad\qquad\;\;\,\qquad\qquad\quad\;\;\, j = 0.
\end{cases}\end{equation}
\end{lemma}
  \begin{proof}
  Looking at the constituent parts of \eqref{chiexpansion}, we start with $(\phi^*\boldsymbol{\omega}^{T})Z\mathbf{f}$. We write
$$(\phi^*\boldsymbol{\omega}^{T})_i = \sum_{j_1 \in \Z} \left(\sum_{k_1 = 0}^{d_x - 1}
    \frac{d_{j_1, k_1}^{(i)}(x)}{r^{j_1}} b_{k_1}^0\right) \frac{dx}{r}.$$ Then
  $\ord_p(d_{j_1, k_1}^{(i)}) \geq \lfloor \frac{j_1}{p} \rfloor + 1$ by \cite[Proof of Proposition~4.9]{Tui17}.
  We have $$f_i = f_{i,0} + f_{i, \infty}  + f_{i, end}\,,$$ where
  $f_{i,0}, f_{i,\infty}$ and $f_{i,end}$ correspond to the three reduction steps (2), (3)
  and (4) in the reduction algorithm from~\cite{Tui17}, summarised
  in~\cite[Algorithm~2.18]{BT19}. By equations (1), (3) and (4) of~\cite{BT19},
there are 
  $\mu_1, \lambda_1\ge 0$ such that
  $$f_{i,0} = \sum_{j_2 = 1}^{\infty}\left(\sum_{k_2 = 0}^{d_x  -1} \frac{c_{j_2,
    k_2}^{(i)}(x)}{r^{j_2}}b_{k_2}^0\right)\,,$$
$$f_{i, \infty} = \sum_{k_3 = 0}^{d_x-1} \sum_{l = 0}^{\mu_1} e_{k_3,
    l}^{(i)} x^l b_{k_3}^0\,, \qquad 
f_{i, end} = \sum_{k_4 = 0}^{d_x-1}\sum_{m=0}^{\lambda_1} u_{k_4, m}^{(i)}x^m b_{k_4}^0\,.$$
    Equation (2) of \cite{BT19} implies the lower bound
  $
  \ord_p(c_{j_2, k_2}^{(i)}) \geq  \left\lfloor \frac{j_2}{p} \right\rfloor + 1 - \log_p
    \lfloor j_2 e_0\rfloor.$
Let 
    \begin{equation}\label{}
    \kappa^{(i)} \colonequals \min(\{0,\ord_p(e_{k_3,l}^{(i)})\} \cup
      \{\ord_p(u_{k_4,m}^{(i)})\})\,\quad \text{ and }\;
  \kappa_1 \colonequals \min_i\{\kappa^{(i)}\}. 
    \end{equation}
  Without loss of generality, the matrix $Z$ has $p$-integral entries. Hence every $(Z\mathbf{f})_i$ is of the form 
\begin{equation}\label{s3}(Z\mathbf{f})_i = \sum_{j_2 = 0}^{\infty} \sum_{k_2 = 0}^{d_x
-1} \frac{g^{(i)}_{j_2,k_2}(x)}{r^{j_2}}b_{k_2}^0\end{equation}
  where for all $k_2$, we have 
    \begin{equation}\label{s31}
      \ord_p(g^{(i)}_{j_2,k_2}) \geq \begin{cases}
      \left\lfloor\frac{j_2}{p}\right\rfloor + 1 - \log_p\lfloor j_2 e_0\rfloor\;,\text{ if }\;j_2>0\\
    \kappa_1\;,\text{ if }\;j_2=0\,.\\
      \end{cases}
    \end{equation}

Let us now consider, for each $i$,  \begin{align*}
  \left(\phi^*\boldsymbol{\omega}^{T}\right)_i \left(Z\mathbf{f}\right)_i  &= \sum_{j_1 \in \Z}\left(\sum_{k_1 =
  0}^{d_x-1} \frac{d_{j_1,k_1}^{(i)}}{r^{j_1}}b_{k_1}^0\right)\left(\sum_{j_2 = 0}^{\infty}\sum_{k_2
  = 0}^{d_x-1}\frac{g^{(i)}_{j_2,k_2}}{r^{j_2}}b_{k_2}^0\right)\frac{dx}{r}\\
&=\sum_{j  = j_1 + j_2 \in \Z,\, j_1 \in \Z,\, j_2 \geq 0} \frac{1}{r^j}\left(\sum_{k =
  k_1 + k_2,\, k_i \in \{0,\ldots,d_x-1\}} \left(d^{(i)}_{j_1,k_1} g^{(i)}_{j_2,k_2}\right)b_k^0\right)\frac{dx}{r}\\
  &\equalscolon\sum_{j\in \Z}\left(\frac{1}{r^j}\sum^{d_x-1}_{k=1} \tau_{jk}b_k^0\right) \frac{dx}{r}.
\end{align*}

We distinguish two cases:
If $j_2 > 0$ then 
    \begin{equation}\label{j2gt0}
   \ord_p(d^{(i)}_{j_1,k_1}g^{(i)}_{j_2,k_2}) \geq
    \left\lfloor\frac{j_1}{p}\right\rfloor + 1 + \left\lfloor\frac{j_2}{p}\right\rfloor +
    1 - \log_p(j_2e_0)
\geq \left\lfloor \frac{j}{p} \right\rfloor + 1 - \log_p((j-1)e_0).
    \end{equation}
If $j_2 = 0$, then $\ord_p(d^{(i)}_{j_1,k_1}g^{(i)}_{j_2,k_2})\geq \left\lfloor
    \frac{j_1}{p}\right\rfloor + 1 + \kappa_1$.
Together, we obtain 
    \begin{equation}\label{j2:0}
\ord_p(\tau_{jk}) \geq \left\lfloor \frac{j}{p}\right\rfloor + 1 - \log_p((j-1)e_0) + \kappa_1.
    \end{equation}

The next term to consider in \eqref{chiexpansion} is
    $\phi^* \eta - p \eta$, where $\eta$ is constructed in~\cite[\S4]{BDMTV19}. Let
    $\kappa_2$ denote the $p$-adic valuation of the vector of coefficients of $\eta$ in
    terms of the basis differentials $\omega_{2g},\ldots,\omega_{2g+2-d}$
    (see~\cite[\S4.1]{BDMTV19}).
Write $$\phi^*\eta - p\eta = \sum_{j \in \Z}\left(\sum_{k = 0}^{d_x-1} \frac{s_{jk}(x)}{r^j}
b_k^0\right) \frac{dx}{r}.$$ Then the $s_{jk}$ satisfy
$\ord_p(s_{jk}) \geq \kappa_2 + \left\lfloor \frac{j}{p} \right\rfloor + 1$ if 
$j\neq 0$ and
$\ord_p(s_{0k}) \geq \kappa_2 +1$, so
\begin{equation}\label{eta-bd}
  \ord_p(s_{jk}) \geq \kappa_2 + \left\lfloor \frac{j}{p}\right\rfloor + 1\;\text{for
  all}\; j\,.
\end{equation}

For the final summand  $\mathbf{g}_0 ^{T}\boldsymbol{\omega}$ in \eqref{chiexpansion} note
that since $F$ has $p$-integral entries, every $(F^TZ \mathbf{f})_i$ has an expansion as in
\eqref{s3}. Because $\omega_i$ is integral for all $i$, the lower bounds in \eqref{s31}
remain valid for $\mathbf{g}_0 ^{T}\boldsymbol{\omega}$. 
The proof of Lemma~\ref{L:chi-bds} follows from this and from~\eqref{j2gt0} and~\eqref{j2:0} upon setting
$\kappa = \min\{\kappa_1,\kappa_2\}$.
  \end{proof}

%%%%%%%%%%%%%%%%%%%%%%%%%%%%%%%%%%%%%%%%%%%%%%%%%%%%%%%%%%%%%%%%%%%%%%%%
We now estimate the precision loss that can occur during the application of
the reduction algorithm from~\cite{Tui17} to the differential $\xi$. 
Our proof is similar to 
the proof of \cite[Prop 4.9]{Tui17}, which estimates the precision loss in the
  reduction of $F^*(\omega_i)$. 
Suppose that $\xi$ is correct to $n$ digits of $p$-adic precision.
First consider terms in~\eqref{chi-sum} with $j>0$. It follows
  from~\eqref{ord-omjk} that  
$j - p\log_p(je_0) \leq pm - p\kappa$
  (note that $\kappa\le 0$).
By~\cite[Prop 3.7]{Tui17}, the
  precision loss at pole order $j$ during the reduction at finite points is at most
  $ \lfloor \log_p(j_{\max}e_0)\rfloor\,,$
  where $j_{\max}$ is the largest integer $j$ such that 
$j - p\log_p(je_0) \leq pn - p\kappa\,.$
  As in the proof of \cite[Prop 4.9]{Tui17}, this might introduce small poles above
  $\infty$, but by the same reasoning as in op. cit., the reduction of these poles leads
  to a loss of precision bounded by $\lfloor\log_p(-(\ord_\infty
  W^{-1})+1)e_{\infty}\rfloor$.
  We set $$g_1(n) \colonequals   
 \lfloor \log_p(j_{\max}e_0)\rfloor + 
\lfloor\log_p(-(\ord_\infty
  W^{-1})+1)e_{\infty}\rfloor\,.$$

If we write 
 $$\xi = \left(\sum_{i=0}^{d_x-1}
 \alpha_i(x,x^{-1})b_i^{\infty}\right)\frac{dx}{r}\quad\text{and}\;\;m_\infty= -\min_i
 \{\ord_\infty\alpha_i-\deg(r) + 1\}\,,$$ 
 then the loss of precision during the reductions above infinity (where $j\leq 0$) is bounded by $g_2\colonequals\lfloor\log_p(m_\infty e_\infty)\rfloor\,.$

Hence we have shown the following:
\begin{lemma}\label{L:}
  Suppose that $\xi$ is correct to $n$ digits of precision. Then
  $\mathbf{c}$ and $\h$ are correct
  to $n-\max\{g_1(n), g_2\}$ digits of precision.
\end{lemma}

\begin{proof}[Proof of Proposition~\ref{P:G}]
  Similar to the $f_i$, we may decompose $\h$ as $\h = \h_0 +
  \h_{\infty} + \h_{end}$, corresponding 
to steps (2), (3) and (4), respectively, in~\cite[Algorithm~2.18]{BT19}.
By the above, the reduction above finite points introduces a denominator of valuation at
  most $\log_p(je_0)$ for pole order $j$, therefore we have
\begin{equation}\label{h-no-loss}
  \h_0 = \sum_{j \geq 1}\sum_{k = 0}^{d_x - 1} \frac{c_{jk}(x)}{r^j} b_k^0\,,\quad\text{ where
  }\;
\ord_p(c_{jk}) \geq \left\lfloor \frac{j}{p} \right\rfloor - 2\log_p(je_0) + \kappa\,.
\end{equation}

Recall that the matrix $G$ is defined in~\eqref{eqn:F-structure}.
 There is no loss of precision when
  evaluating $\mathbf{f}(P)$  by~\cite[Prop. 4.5]{BT19}.
  By our assumption that $F$ and $Z$ are $p$-integral, there is no precision loss when evaluating $\mathbf{g_0}(P)$.
  Using the bounds~\eqref{h-no-loss}, the proof of~\cite[Prop. 4.5]{BT19} shows that
  $\h(P)$ is accurate to $n$ digits of precision as well.
  Since $\mathbf{g} = \mathbf{g_0} + \mathbf{c}$, the proposition
  follows.
\end{proof}

%%%%%%%%%%%%%%%%%%%%%%%%%%%%%%%%%%%%%%%%%%%%%%%%%%%%%%%%%%%%%%%%%%%%%%%%
\subsubsection{Frobenius-equivariant splitting for general points}\label{subsec:gen}

For $x\in X(\Q_p)\cap ]\mathcal{U}[$, not necessarily Teichm\"uller,
the Frobenius-equivariant splitting $\lambda^\phi(x)$ of $M(x)$ is given by 
\begin{small}
  \begin{equation}\label{frobsplit}
\begingroup
\setlength\arraycolsep{4.5pt} 
\begin{pmatrix}
1 & 0 & 0 \\
\int_{x}^{x_0} \boldsymbol{\omega} & 1 & 0 \\
\int_{x}^{x_0} \eta+\int ^x _{x_0 }\boldsymbol{\omega }^{\intercal }Z \boldsymbol{\omega } & \int_{x}^{x_0} \boldsymbol{\omega}^{\intercal}Z & 1 \\
\end{pmatrix}\cdot
\begin{pmatrix}
1 & 0 & 0 \\
\int_{b_0}^b \boldsymbol{\omega} & 1 & 0 \\
\int_{b_0}^b \eta +\int ^b _{b_0 }\boldsymbol{\omega }^{\intercal }Z \boldsymbol{\omega } & -\int_{b_0}^b \boldsymbol{\omega}^{\intercal}Z & 1 \\
\end{pmatrix}\cdot
  \lambda^\phi(x_0),
\endgroup
\end{equation}
\end{small}\\
where $x_0$ is the Teichm\"uller point in the disk of $x$.
The first two matrices in~\eqref{frobsplit} correspond to parallel transport of $\Lambda$ from $x$ to
$x_0$ and from $b_0$ to $b$, respectively.

%%%%%%%%%%%%%%%%%%%%%%%%%%%%%%%%%%%%%%%%%%%%%%%%%%%%%%%%%%%%%%%%%%%%%%%%

For the local height $h_p(A(x))$, we need the Frobenius-equivariant
splitting $\lambda^\phi(x)$ both for fixed $x$
and for $x$ varying inside a residue disk.
We start by bounding the valuations of the coefficients of power series expansions of the
differentials in the parallel transport matrices of $\Lambda$ in terms of a local
coordinate $t$ at a fixed affine point $y_0 \in X(\Q_p)\cap ]\mathcal{U}[$.
By assumption, the entries of the
expansions of 
$\boldsymbol{\omega}$ and $\boldsymbol{\omega}^{\intercal}Z$ all have
integral coefficients, so their integrals have entries whose $i$-th
coefficient has valuation $\ge -\ord_p(i)$. 
Therefore, we have
\begin{equation}\label{}
  \boldsymbol{\omega}(t)^{\intercal}Z\int \boldsymbol{\omega}(t) =
  \sum_{i\ge 1} a_it^i\,,\quad\text{where }\;\ord_p(a_i)\ge
  -\lfloor \log_p(i)\rfloor\,.
\end{equation}
It follows that 
\begin{equation}\label{double-int-coeffs}
  \int(\boldsymbol{\omega}(t)^{\intercal}Z\int \boldsymbol{\omega}(t)) =
  \sum_{i\ge 1} b_it^i\,,\quad\text{where }\;\ord_p(b_i)\ge
  -2\lfloor \log_p(i)\rfloor\,.
\end{equation}

By construction, the coefficients of $\eta$ in terms of
$\omega_{2g},\ldots,\omega_{2g+d-2}$ are polynomials in $x$.
Define $d_i(\eta)$  to be the 
valuation of the $i$th coefficient if $i$ is smaller than the maximum of the  degrees of
these coefficients and~0 otherwise.
Then the $i$th coefficient of the integral of $\eta$ has valuation $\ge
-\ord_p(i)-d_i(\eta)$.
Hence, the $i$th coefficient of every expansion of the parallel transport matrix in $t$
has valuation at least
\begin{equation}\label{bdpt}
  \varphi(i) \colonequals -\lfloor \log_p(i)\rfloor+ \min\{d_i(\eta), -\lfloor \log_p(i)\rfloor\}\,.
\end{equation}

For definite parallel transport from $y_0$ to another $\Q_p$-point $y_1$ in the same residue disk,
we need to evaluate the integrals above. Suppose that $y_0,y_1$, and the coefficients of the
expansions of $\boldsymbol{\omega}$ and $\eta$ are correct to $n$ digits of $p$-adic
precision, and suppose that the expansions are truncated modulo $t^l$.
Let 
\begin{align*}
  \nu_1\colonequals 1+\min_{i\ge l}\{i-\lfloor\log_p(i+1)\rfloor\}\quad\text{and}\;
  \nu_2\colonequals n+\min_{0\le i\le l-1}\{i-\lfloor\log_p(i+1)\rfloor\}\,.
\end{align*}
Then 
$\int_{y_0}^{y_1}\omega_j$ and
$\int_{y_0}^{y_1}(Z\boldsymbol{\omega})_j$ are correct to $\min\{\nu_1,\nu_2\}$ digits by~\cite[Prop.~4.1]{BT19}. 
The proof of~\cite[Prop.~4.1]{BT19} requires that the differential we integrate has integral
coefficients. A modification of this proof yields that the integral $\int_{y_0}^{y_1}\eta$ is correct to 
$
\min \{\nu_1', \nu_2\}
$
digits, where
$\nu'_1=1+\min_{i\ge l}\{i-\lfloor\log_p(i+1)\rfloor-d_i(\eta)\}.$
A similar modification shows that the double integral 
$\int ^{y_1} _{{y_0}}\boldsymbol{\omega }^{\intercal }Z \boldsymbol{\omega } $ is correct to
$\min\{\nu_1'', \nu_2'\}$ digits, where
$
  \nu_1''= 1+\min_{i\ge l}\{i-2\lfloor\log_p(i+1)\rfloor\}$ and
 $$ \nu_2'= n-\lfloor\log_p(n)\rfloor+\min_{0\le i\le
 l-1}\{1-\lfloor\log_p(i+1)\rfloor\}\,.$$ 
Hence we obtain the following:
\begin{lemma}\label{L:pt-eval}
  The parallel transport matrix from ${y_0}$ to ${y_1}$ is correct to 
$
  \min\{\nu_1',\nu''_2, \nu_2'\}
$
digits of precision.
\end{lemma}

Using \eqref{frobsplit}, we can finally bound the loss of precision
in the computation of $\lambda^\phi(x)$ for fixed points
$x\in X(\Q_p)\cap ]\mathcal{U}[$ by combining Lemma~\ref{L:pt-eval} and Proposition~\ref{P:G}.

%%%%%%%%%%%%%%%%%%%%%%%%%%%%%%%%%%%%%%%%%%%%%%%%%%%%%%%%%%%%%%%%%%%%%%%%

\subsection{Global heights}\label{subsec:globht}
We now discuss the possible precision loss in the computation of the
global height $h$.
In Step~\eqref{solve} of Algorithm~\ref{alg:qc} we solve for $d_1,\ldots,d_g$ such that 
\begin{equation}\label{htssolver}
h  =\sum_{i}d_i\Psi_i\
\end{equation}
in terms of a basis $\{\Psi_i\}$ of bilinear pairings on
$\HH^0(X_{\Q_p}, \Omega^1)^{\vee}$ 
 by evaluating $h$ and the $\Psi_i$.
Recall that our method for determining the coefficients depends on whether there are
sufficiently many rational points on $X$ in the sense of~\S\ref{subsec:global_height}.
If this is the case, meaning that we can use a basis of consisting of 
$\pi(\mathrm{A}_Z(b,x))$ for rational points $x\in X(\Q)\cap ]\mathcal{U}[$,
then we need to compute $h_p(\mathrm{A}_Z(b,z)$ and
$\pi(\mathrm{A}_Z(b,x))$, and then apply simple linear algebra. The
precision loss in the computation $h_p(\mathrm{A}_Z(b,x))$ has already been
bounded and $\pi(\mathrm{A}_Z(b,x))$ can be obtained directly from
the same data (see~\cite[Equation~(41)]{BBBLMTV19}).
The loss of precision in the linear algebra computations 
is easy to detect in practice, so we do not bound it explicitly here.

In the other case, the basis $\Psi_i$ is given in terms of products of
abelian integrals. As mentioned above, the loss of precision in their
computation is estimated in~\cite{BT19}. It remains to discuss precision loss in the computation of Coleman--Gross local heights $h_p(D_1,D_2)$, where $D_1, D_2$ are divisors in $\Div^0(X)(\Q_p)$ for $X$ a hyperelliptic curve subject to the hypotheses of Algorithm \ref{alg:hp}; see \cite[\S 6.2]{bb:heights} for further details. Choosing $\omega$ in Step (1) can be done up to the precision of the points in the support of the divisor $D_1.$  To compute $\Psi(\omega)$ and $\omega_{D_1}$ to $O(p^n)$ in Step (2),  see Section 5.2 and Section 6.2.3 of \cite{bb:heights}: one needs to compute the local coordinates $(x(t),y(t))$ at infinity, with $x(t)$ to precision $t^{2(2g-1)}$ and $y(t)$ to precision $t^{2g-1}$, where these $t$-adic estimates are made based on the maximal pole order in the basis of $H^1_{\dR}(X)$. Step (4) proceeds similarly to this step as well. 

In Step (5), the tiny integrals are computed as in \cite[\S 6]{bb:heights}. In previous steps, we wrote $\Psi(\alpha)$ and $\Psi(\beta)$ as $\Q_p$-linear combinations of the basis elements of $H_{\dR}^1(X)$, up to precision $O(p^n)$. Note that the hypothesis that $p> 2g-1$ is to ensure that the cup product matrix has entries that are $p$-integral, so no precision loss comes from the cup product matrix.  

Finally, for $\sum_{A \in X(\C_p)}\Res_{A}\left(\alpha\int\beta\right)$, we consider the cases of $A$ a non-Weierstrass point (where we describe the computation in the annulus of $A$) versus $A$ Weierstrass (where we have just one contribution, at the Weierstrass point).  
If $A \neq (0,0)$ is a Weierstrass point, we compute the local coordinate
$(x(t),y(t))$ at $A$ to precision $t^{2pn-p-1}$ (see the corrected
Proposition 6.5 in \cite{bb:heightserrata}) so that
$\Res_A(\alpha\int\beta)$ is computed to $n$ digits of $p$-adic
precision. 

Now we consider the non-Weierstrass poles of $\alpha$.  For the
annulus of a non-Weierstrass pole $A$, the generic situation is handled by \cite[Corollary 6.4]{bb:heights}.
By \cite[Remark~4.10]{bb:heights}, we consider all $A \in \{P_i,
Q_j\}_{i,j}$
where $x(P_i)$ corresponds to a root of an irreducible factor of $x^p -
x(P)$ (and similarly where $x(Q_j)$ corresponds to a root of an irreducible
factor of $x^p - x(Q)$). For these $i,j$, we compute $\int_P^{P_i}\beta$ and  $\int_Q^{Q_j}
\beta$ and trace down to $\Q_p$.
We suppose $P \in X(\Q_p)$ has  precision $O(p^n)$. Fix $m$ and
suppose $\beta$ is computed to $t^{d_i m}$ at $P_i$, where $d_i =
[\Q_p(P_i): \Q_p]$. Let $\pi_i$ be a uniformiser of
$\Q_p(P_i)$. Note that $P$ is known to $d_i n$ $\pi_i$-adic digits, and
suppose that $P_i$ is known to $n_i$ $\pi_i$-adic digits. Then the
$\pi_i$-adic precision of $\int_P^{P_i} \beta$ is at least $\min\{n_i, d_i
n, \lfloor d_i m + 1\rfloor - \log_p(d_i m + 1)\}.$ We similarly repeat
this for $Q$ and the corresponding $Q_j$. Hence
$\sum_A\Res_A(\alpha\int\beta)$, where the sum is over all non-Weierstrass
poles $A$ of $\alpha$, is correct to $p$-adic precision 
\[
\min_{i,j} \{n_i, d_i n, \lfloor d_i m + 1\rfloor - \log_p(d_i m + 1),
  n_j, d_j n, \lfloor d_j m + 1\rfloor - \log_p(d_j m + 1)\}\,, 
\]
where we consider the corresponding $i,j$ for all $P_i$ and all $Q_j$.

%%%%%%%%%%%%%%%%%%%%%%%%%%%%%%%%%%%%%%%%%%%%%%%%%%%%%%%%%%%%%%%%%%%%%%%%
\subsection{Coefficients of the quadratic Chabauty function and root finding}\label{subsec:boundcoeffs}
The previous results of this section bound the loss of precision in
the computation of the quadratic Chabauty function $\rho = h-h_p$. Let $D\subset X(\Q_p)\cap]\mathcal{U}[$ be a residue disk and let $x_0$ be the Teichm\"uller point in $D$.
We now bound the valuations of the coefficients of the expansion of
$\rho$ in $D$ and show how to provably compute its roots to desired
precision.

In our algorithm, we fix a point $x_1\in D$, and we compute the
Frobenius-equivariant splitting $\lambda^\phi(x)$  on $D$ as a power
series in a local coordinate $t$ in $x_1$ by first computing
$\lambda^\phi(x_1)$ from $\lambda^\phi(x_0)$ and then multiplying this by the parallel
transport matrix from $x_1$ to $x$.
To bound the valuations of the coefficients of the entries of
$\lambda^\phi(x)$,

we first compute 
\[
c_1\colonequals
\ord_p(\lambda^\phi(x_1)) 
\]
using Lemma~\ref{L:pt-eval}.
By the above, we find that the $i$th coefficient of every entry of the expansion of $\lambda^\phi(x)$ has valuation at least 
$ \varphi(i) + c_1$.  
We use this to bound the valuations of the coefficients of the local height $h_p$. 
Recall from~\S\ref{subsec:global_height} that we use a height with respect to an
$\End(J)$-equivariant splitting of the Hodge filtration; let
$v_{\spl}$ be the smallest valuation of the coefficients of this splitting in terms of our basis
$\boldsymbol{\omega }$. 
We denote by $\ord_p(\gamma_{\Fil})$ the smallest valuation in the
coefficients of the rational function $\gamma_{\Fil}$, and we set
$$
c_2 \colonequals \min\{0,v_{\spl}, \ord_p(\beta_{\Fil}),v_{\spl} +\ord_p(\beta_{\fil}) \}\,.
$$

\begin{lemma}\label{L:hpbd}
  Let 
  $$h_p(x(t)) = \sum_{i\ge 0} h_it^i$$ be
  the expansion of $h_p$ on the residue disk $D$ in the local parameter
  $t$. 
  Then we have
\begin{equation}\label{}
  \ord_p(h_i) \ge \min\{\ord_p(\gamma_{\Fil}), \varphi(i) +c_2\}.
\end{equation}
\end{lemma}
\begin{proof}
  This follows from the discussion above and from~\eqref{hpformula}, 
    which expresses $h_p(x)$ in terms of $\lambda^{\Fil}(x)$ and
    $\lambda^\phi(x)$.
\end{proof}

We set $c_3\colonequals  \min_i\{\ord_p(d_i)\}$, where the $d_i$ are the coefficients
in~\eqref{htssolver}.
Let $i_0\ge0$ be such that
$$-\lfloor\log_p(i)\rfloor \le \min\left\{d_{i}(\eta), \left\lfloor
\frac{\ord_p(\beta_{\Fil})}{2}\right\rfloor, \left\lfloor
\frac{\ord_p(\gamma_{\Fil})-c_2}{2}\right\rfloor \right\}
$$
for all $i\ge i_0$.
Then we have $\varphi(i) = -2\lfloor \log_p(i)\rfloor+c_1$ for all $i\ge i_0$. This
proves the following:
\begin{prop}\label{P:rhobd}
  Let 
  $$\rho(t)  =\sum_{i\ge 0} \rho_it^i $$
be the expansion of the quadratic Chabauty function $\rho=h-h_p$ on $D$.
  If $i\ge i_0$, then we have
$$
  \ord_p(\rho_i) \ge -2\lfloor\log_p(i)\rfloor + c_1+ \min\{c_2, c_3\}.
$$
\end{prop}
Together with Proposition~\ref{P:rhobd}, the following result allows us to provably determine the roots of
$\rho$ to any desired precision.
\begin{lemma}
Suppose $F(x)=\sum _{i\geq 0}F_i x^i \in \Q _p [\! [x]\! ]$ is such that there are integers $k,m,n$ satisfying  
\[
  \min \{\ord_p (F_i )+i:i\geq 0\} =k
\]
and 
\[
  \max \{i\geq 0:\ord_p (F_i )+i=n \} <m,
\]
  and furthermore that $F$ has at most $d$ roots in the closed disk $\{
    \ord_p (x)\geq 1\}$. 
  Then the roots of $F $in the ball $\{\ord_p (x)\geq 1\}$ can be determined, with multiplicity, to precision $(n-k)/d$, by computing $F_0 ,\ldots ,F_{m-1} $ modulo $p^n$.
\end{lemma} 
\begin{proof}
By our assumptions, $F(px)$ lies in $p^k \mathbf{Z}_p [\! [x]\! ]-p^{k+1}\mathbf{Z}_p [\!
  [x]\! ]$. Hence the power series $G(x)\colonequals p^{-k}F(px)$ lies in $\mathbf{Z}_p [\! [x]\! ]-p\mathbf{Z}_p [\! [x]\! ]$. Furthermore, by our assumptions, for any $\alpha \in \mathbf{Z}_p $, the positive slopes of the Newton polygon of $G(x+\alpha )$ are uniquely determined by the first $m$ coefficients. If $G(x)$ is congruent modulo $p^{n-k}$ to a polynomial $H$ in $\mathbf{Z}_p [x]$ with a root $\alpha \in \overline{\mathbf{Z}}_p $ of multiplicity $e$, then the valuation of the first $e$ coefficients of $G(x+\alpha )$ must be at least $n-k$. Since $G(x+\alpha )$ has degree $\leq m$ mod $p^{n-k}$ and has at least one coefficient of valuation zero, we deduce that the Newton polygon of $G(x+\alpha )$ must contain a segment of slope $\geq (n-k)/d$ of length at least $e$.
\end{proof}
\begin{Remark}
In practice, we usually apply this with $d=1$, by recentering and rescaling
  our power series so that there is only one root in the ball $\{
    \ord_p (x) \geq 1\}$ (and because in practice the power series do not typically have repeated roots). Hence most loss of precision occurs from $k$ being large, rather than $d$.
\end{Remark}

%=========================================================================
\section{Examples}\label{sec:examples}

%=========================================================================
In this section, we apply our techniques to compute the rational points on
\begin{itemize}
  \item  the exceptional modular curve $X_{S_4}(13)$ (see \S\ref{subsec:s4});
  \item all curves $X_0^+(N)$ of genus 2 and 3 for which $N$ is prime and the rational
    points were not previously known (see \S\ref{subsec:x0n+}); 
  \item two genus~2 curves of interest in Mazur's Program B (see \S\ref{subsec:dzb});
  \item two genus~2 curves with Jacobian of $\mathrm{GL}_2$-type that have nontrivial local
    height contributions away from $p$ (see \S\ref{subsec:nontriv_away});
    \item the non-split Cartan curve $X_{\ns}^+(17)$ (see \S\ref{sec:ns17}).
\end{itemize}
For the computations, we used our {\tt Magma} implementation.
The code used for the examples, along with log files, can 
be found in the folder \verb|Examples| at~\cite{QCCode}.

%%%%%%%%%%%%%%%%%%%%%%%%%%%%%%%%%%%%%%%%%%%%%%%%%%%%%%%%%%%%%%%%%%%%%%%%
\subsection{The exceptional curve $X_{S_4}(13)$}\label{subsec:s4}
Recall that for a prime $\ell \geq 5$, any proper subgroup of $\mathrm{GL}_2(\F_{\ell})$ is conjugate to a subgroup of a Borel subgroup, the normaliser of a Cartan subgroup, or an ``exceptional'' subgroup with projective image isomorphic to $S_4, A_4$, or $A_5$. The field of definition of the modular curves attached to the exceptional subgroups is the unique quadratic subfield $\Q(\sqrt{\pm \ell})$ of the cyclotomic field $\Q(\zeta_{\ell})$, with the exception of the curves $X_{S_4}(\ell)$ for $\ell \equiv \pm 3 \pmod{8}$, which are defined over $\Q$. For such values of $\ell$, we would therefore like to determine $X_{S_4}(\ell)(\Q)$. 

\par Serre~\cite{Ser72} shows by a monodromy argument that such tetrahedral modular curves have no points defined over $\Q_{\ell}$ when $\ell$ is large enough, and in particular he obtains 
\[
X_{S_4}(\ell)(\Q) = \emptyset, \qquad \mbox{if} \ \ \ell > 13.
\]
The curves $X_{S_4}(3)$ and $X_{S_4}(5)$ are both of genus zero, and contain a unique rational cusp. Ligozat \cite{Lig76} showed that $X_{S_4}(11)$ is an elliptic curve of conductor $11^2$ whose Mordell--Weil group is trivial, where the unique rational point is CM, corresponding to discriminant $D=-3$. This leaves only the curve $X_{S_4}(13)$, which has genus $3$. In fact, this curve is the last remaining modular curve of level $13^n$ whose rational points have not been determined. 

Using modular symbols algorithms, Banwait--Cremona \cite{banwait-cremona} show that the curve $X_{S_4}(13)$ is a smooth plane quartic whose canonical model is given by
\begin{align*}
4x^3y - 3x^2y^2 + 3xy^3 - x^3z + 16x^2yz - 11xy^2z + \\ 5y^3z + 3x^2z^2 + 9xyz^2 + y^2z^2 + xz^3 + 2yz^3 = 0.
\end{align*}

Furthermore, they exhibit the following four rational points
\[
\{  \left(1 : 3 : -2\right), \left(0 : 0 : 1\right), \left(0 : 1 :
0\right),\left(1 : 0 : 0\right)\} \subseteq X_{S_4}(13)(\Q)\,,
\]
where the rational point $(0 : 0 : 1)$ corresponds to an elliptic curve with CM by the order of discriminant $D = -3$, and the three other rational points correspond to non-CM elliptic curves over $\Q$ with projective mod $13$ image equal to $S_4$, whose $j$-invariants are given by
\[
\begin{array}{lll}
j & = &  \displaystyle \phantom{-\, } \frac{2^4\cdot 5 \cdot 13^4 \cdot 17^3}{3^{13}} \qquad j \ = \ \displaystyle - \, \frac{2^{12}\cdot 5^3 \cdot 11 \cdot 13^4}{3^{13}} \\ [12pt]
j & = & \displaystyle \phantom{-\, }  \frac{2^{18}\cdot 3^3 \cdot 13^4\cdot
  127^3 \cdot 139^3 \cdot 157^3 \cdot 283^3 \cdot 929}{5^{13} \cdot
  61^{13}}\,. \\
\end{array}
\]

The Jacobian of $X_{S_4}(13)$ is isogenous to that of $X_{\rm s}^+(13)$, so it is absolutely simple and has Mordell--Weil rank~3 over $\Q$ by the results of~\cite[\S 6]{BDMTV19}. The curve has potential good reduction at $p = 13$, as can be seen, for instance, using the {\tt Sage} toolbox {\tt MCLF}. 

We determine the set of rational points on the curve $X_{S_4}(13)$ using quadratic
Chabauty with $p=11$ for the affine patches
\begin{align*}
  y^4 &+ (18x + 9)y^3 + (160x^2 + 176x + 52)y^2 + (560x^3 + 832x^2 + 384x +
    48)y\\& + 192x^4 + 512x^3 + 384x^2 + 64x = 0
\end{align*}
    and 
    \begin{align*}
      y^4& + (9x + 9)y^3 + (52x^2 + 72x + 36)y^2 + (48x^3 + 240x^2 + 208x + 
      64)y \\&+ 64x^3 + 192x^2 - 64x = 0.
    \end{align*}
The computation is analogous 
to the computation of $X^+_{\rm s}(13)(\Q)$ in~\cite{BDMTV19}. 
The Hecke operator $T_{11}$ generates the
Hecke algebra, as can be verified, for instance, by checking the analogous statement for $X^+_{\rm s}(13)$. Hence we may construct suitable cycles $Z_1, Z_2$ from
$T_{11}$ and its square, respectively. 
The set of common zeroes of the resulting
quadratic Chabauty functions consists precisely of the known rational points, so we
obtain Theorem~\ref{T:XS4}. 

In order to solve for the height pairing, we use the~4
known rational points and the cycle $Z_1$, so the resulting function
automatically vanishes there. However, since the cycles $Z_1$ and $Z_2$ are
independent, and $Z_2$ is not used to solve for the height, the vanishing of the resulting
function in the rational points provides a check for the correctness of our code. 

\begin{Remark}
Since the Jacobian of $X_{S_4}(13)$ is isogenous to that of $X^+_{\mathrm{s}}(13)$, even if there were not enough rational points on $X_{S_4}(13)$ to solve for the height pairing, one could instead solve for it using $X^+_{\mathrm{s}}(13)$.
\end{Remark}

\subsection{The Atkin--Lehner quotients $X_0^+(N)$ } \label{subsec:x0n+}

\par For  a positive integer $N$, consider the Atkin--Lehner involution $w_N$ acting on
the modular curve $X_0(N)$. Then the 
quotient 
\[
X_0^+(N ) \colonequals X_0(N)/\langle w_N \rangle
\]
is a smooth projective curve defined over $\Q$ whose non-cuspidal points classify unordered pairs $\{E_1,E_2\}$ of elliptic curves admitting an
$N$-isogeny between them. 
The study of rational points on these curves is also important in an
ongoing research program aiming to compute quadratic points on the modular curves
$X_0(N)$; see, for instance, recent work of Box~\cite{Box19}. 
Among the rational points, we distinguish between cusps, CM-points and {\em exceptional
points}, those which are neither cusps nor CM points. The exceptional points correspond to quadratic $\Q$-curves without CM.

In this section, we restrict to prime values $N$ such that $X_0^+(N)$ has genus 2 or 3. 
Galbraith~\cite{Gal96} has computed models for all these curves (and many more) by finding
relations in the vector space spanned by the newforms of level $N$ and weight $2$ that are
invariant under $w_N$. Up to conjugation, there is a unique such newform. 

By work of Ogg, for prime level $N$, the curve $X_0^+(N)$ has genus 2 if and only if
\begin{equation}\label{eq:X0N+g2}
  N \in \{67, 73, 103, 107, 167, 191 \}\,.
\end{equation}
It has genus 3 if and only if
\begin{equation}\label{eq:X0N+g3}
    N \in \{97, 109, 113, 127, 139, 149, 151, 179, 239\}\,.
\end{equation}
Models for all these curves were communicated to us by Elkies; 
one can also find such models in Galbraith's
thesis~\cite{Gal96} or by
using the Magma command \verb|X0NQuotient|.

Via a search for small rational points, Galbraith~\cite{Gal99} found exceptional rational
points on $X_0^+(N)$ for $N=73, 91, 103, 191$ (genus~2) and $N=137, 311$ (genus~4). The latter
examples disproved an earlier conjecture of Elkies that 
there are no exceptional rational points on non-hyperelliptic $X_0^+(N)$ for prime level $N$. 
In~\cite{Gal02}, Galbraith also finds an exceptional point on $X_0^+(125)$ 
and conjectures that there are no further exceptional points on modular curves
$X^+_0(N)$ of genus $2\le g\le 5$.

Together with~\cite{BBBLMTV19} and~\cite{Gal96}, our computations described below prove
Theorem~\ref{T:X0N+}. 
We first check that for level $N$ as in~\eqref{eq:X0N+g2} and~\eqref{eq:X0N+g3} the curves
$X_0^+(N)$ satisfy
the requirements to apply our algorithm.
The Jacobian $J_0^+(N)$ of $X_0^+(N)$ has real multiplication over $\Q$, so the Picard
number is at least $g$.
Using \verb|Magma| we computed the $L$-function of the corresponding newforms to show that the analytic rank is $g$, so 
the work of Gross--Zagier and Kolyvagin--Logachev proves that the rank of $J(\Q)$ is
exactly $g$.  For the genus~2 examples, we also applied two-descent on
$J_0^+(N)$, as implemented in~{\tt Magma}, to have an independent check. 

The curves $X_0^+(N)$ have good reduction away from $N$, but in contrast to $X^+_{\rm
ns}(13)$ and $X_{S_4}(13)$, they do not have potentially good reduction at $N$.
Nevertheless, the following result implies that when applying quadratic Chabauty, there
are no nontrivial contributions to the height away from $p$.  

\begin{lemma}\label{L:X0N+model}
There is a regular semi-stable model $\mathcal{X}_0^+(N)$ of $X_0^+(N)$ over $\Z_N$ whose special fibre has a unique irreducible component. In particular, the local height $h_N$ is trivial on $X_0^+(\Q_N)$. 
\begin{proof}
If $N=2,3$ the result is readily checked. When $N > 3$ the Atkin-Lehner quotient of the model $\mathcal{X}_0(N)$ for $X_0(N)$ over $\mathrm{Spec} \Z$ constructed by Deligne--Rapoport \cite{DR73} is shown by Xue \cite{Xue09} to be regular and semi-stable. Its special fibre at $N$ is a projective line, with an ordinary double point for every conjugate pair of supersingular $j$-invariants in $\F_{N^2} \backslash \F_{N}$. It follows from Theorem~\ref{thm:BeD} that $h_N$ is trivial. 
\end{proof}
\end{lemma}

Finally, we checked for all $N$ in~\eqref{eq:X0N+g2} and in~\eqref{eq:X0N+g3} that the Jacobian is absolutely simple by finding a prime $q$ of good
reduction such that $J_{ \F_q}$ is absolutely simple, using the criterion of
Howe and Zhu~\cite[Proposition~3]{HZ02}.

%%%%%%%%%%%%%%%%%%%%%%%%%%%%%%%%%%%%%%%%%%%%%%%%%%%%%%%%%%%%%%%%%%%%%%%%
\subsubsection{Genus 2}\label{subsec:x0n+g2}
In~\cite{BBBLMTV19}, the rational points on $X_0^+(N)$ for $N=67,73,103$
were computed.
Using a combination of quadratic Chabauty and the Mordell--Weil sieve, it is shown there
that $X_0^+(67)(\Q)$ contains no exceptional points and that the sets
$X_0^+(73)(\Q)$ and $X_0^+(103)(\Q)$ both contain one pair of exceptional
points each, with respective $j$-invariants (see~\cite[Table~1]{Gal99})
\[
\begin{array}{l}
j= (81450017206599109708140525\pm 14758692270140155157349165\cdot\sqrt{-
127})/ 2^{74}, \\
  j = (35982263935929364331785036841779200 \\\qquad \pm
669908635472124980731701532753920 \cdot\sqrt{5\cdot577}.
\end{array}
\]
The remaining prime level genus 2 curves $X_0^+(107)$, $X_0^+(167)$, and $X_0^+(191)$ are more challenging, because
they do not have sufficiently many rational points in the sense of~\S\ref{subsec:global_height} to solve for the height pairing, so we need to compute heights between divisors. 
In all cases, the quadratic Chabauty function $\rho=h-h_p$ has $p$-adic zeroes that do not come
from a rational point; to verify this, we apply 
the Mordell--Weil sieve.

\begin{Example}\label{E:107}
  We discuss our computations for the example $X\colonequals X_0^+(107)$  in some detail.

  We look for a prime $p$ of good reduction such that 
  \begin{itemize}
    \item there is a unique $\Q_p$-rational Weierstrass disk, and it does not contain 
      known rational points,
    \item the Hecke operator $T_p$ generates the Hecke algebra, and
    \item $p$ is suitable for the Mordell--Weil sieve.
  \end{itemize}

For $p=61$, the first two conditions are satisfied; moreover, we have
  \[
    J(\F_{229})\simeq \faktor{\Z}{(4\cdot 61)\Z}\times \faktor{\Z}{(4\cdot 61)\Z}\,,
  \]
  and since $J(\F_{61})\simeq \Z/(31\cdot 151)\Z$ has quite smooth order,
  $61$ is a 
  suitable prime. We now go through the steps in
  Algorithm~\ref{alg:qc}, applied to $X$.

  Step~\eqref{smallpts}
  The model
  \[
  y^2 = -3 x^6 - 4 x^5 - 2 x^4 + 2 x^3 + 5 x^2 + 2 x + 1 \equalscolon f(x),
\]
  of $X$ has~6 small rational points of exponential height at
  most~1000, given by $\{(0,\pm1), (\pm1, \pm1) \}$. It also has
no $\Q_{61}$-adic points at infinity, so that we only need to run our algorithm
  for one affine patch. We fix the base point $b=(0,-1)$.

  Steps (\ref{intsympbasis}, \ref{frobhecke}) are exactly as in~\cite{BDMTV19}. 
  
  Step \eqref{hodgeequiv}: We may use the unit root splitting, since $p=61$ is ordinary. (See Remark \ref{rem:unitroot}.)

  Step \eqref{corrs}: Using Step~\eqref{frobhecke},  we find for  
    \begin{equation}\label{eq:tr0}
      Z=Z_1= (\Tr(T_{{61}})\cdot I_{4}- 4T_{{61}})C^{-1} =
\left( \begin{array}{cccc}
  0& 2/3   &-2  & 4\\
  -2/3 &   0   &4  &  2\\
   2  &  -4  &  0   & 0\\
   -4  & -2 &   0   & 0\\
\end{array} \right)\,,
    \end{equation} that
$$Z=\sum_{i,j}  Z_{ij} \omega_i\otimes\omega_j\in \HH^1_{\dR}(X/\Q_{61})\otimes
  \HH^1_{\dR}(X/\Q_{61})\,$$ 
  corresponds to a nontrivial cycle $Z \in \ker(\NS(J)\longrightarrow \NS(X))$
  where $C$ is the standard symplectic matrix of dimension $2g$ and $\boldsymbol{\omega}$ is the
  basis found in Step~\eqref{intsympbasis}.

  Step~\eqref{aff}: The Hodge filtration for $Z$ is given by $\gamma_{\Fil}=-4x-4$ and
$\beta_{\Fil}=0$. After computing the Frobenius structure, we obtain a power series
  expansion of the function $x\mapsto h_{61}(\mathrm{A}(x))$ on all residue disks of $X(\Q_{61})$, except
for the disks at infinity and the unique Weierstrass disk containing points that reduce
to $(31,0)$. 

Step~\eqref{solve}: The points $P, Q\in J(\Q)$ with respective Mumford
  representations $(x^2 + x ,1)$ and $(x^2+1, 2x-1)$ generate a
  subgroup of $J(\Q)$ of index~2. To solve for the height pairing via~\S\ref{subsec:coleman-gross}, we need divisor
  representatives with support in distinct non-Weierstrass residue disks.
  Let $E$ be the degree~2 divisor on
  $X$ cut out by the functions $x^2+1$ and $2x-1$ and let $E'$ be its image under the hyperelliptic involution. We set
\begin{align*}
  D_1 &= (0,1) + (-1,1) - \div_0(x-1),\quad 
  D'_1 = (0,-1) + (-1,-1) - \div_0(x-7) \\
  D_2 &= E-\div_0(x-7), \quad
  D'_2 =  E' - \div_0(x-1)\,.
\end{align*}
Then we have
  $h(P,Q) = \sum_v h_v(D_1, D_2)$ and 
\begin{align*}
  h(P,P) &= -\sum_v h_v(D_1, D_1')\,,\quad
h(Q,Q) = -\sum_v h_v(D_2, D_2')\,.
\end{align*}
The divisors above all split over $\Q_{61}$, so we can compute the height pairings
$h_{61}(D_1, D_2)$, $h_{61}(D_1, D'_1)$ and $h_{61}(D_2, D'_2)$, 
 working on a monic odd degree model over $\Q_{61}$.
Using {\tt Magma}'s implementation of the algorithm described
in~\cite{Mul14}, we also find 
\begin{align*}
  \sum_{\ell\ne {61}} h_{\ell}(D_1, D'_1) &= -2\log_{61} 2+2\log_{61} 3-\log_{61} 7\,,\\ 
  \sum_{\ell\ne {61}} h_{\ell}(D_1, D_2) &= 2\log_{61} 2-2\log_{61} 3 + \log_{61} 7,\\ 
  \sum_{\ell\ne {61}} h_{\ell}(D_2, D'_2)& = 3\log_{61} 2-\log_{61} 5\,,
\end{align*}
and we conclude that
$$
  h= \alpha_{00}g_{00} + \alpha_{01}g_{01}+ \alpha_{11}g_{11}\,,$$
  where 
\begin{align*}
  \alpha_{00} &=58\cdot 61^{-1} + 19 + 2\cdot 61 + 43\cdot 61^2 + O(61^3)\,\\
  \alpha_{01} &=  43\cdot 61^{-1} + 48 + 44\cdot 61 + 41\cdot 61^2 + O(61^3)\,\\
  \alpha_{11} &=49\cdot 61^{-1} + 13 + 55\cdot 61 + 2\cdot 61^2 +
  O(61^3)\,,
\end{align*}
and the $g_{ij}$ are defined in~\eqref{E:gij}.

Steps~\eqref{roots} -- \eqref{return}: Combining the functions resulting from
Steps~\eqref{aff} and~\eqref{solve}, we find a power series expansion
of the quadratic Chabauty function
$$\rho = h - h_{61}\colon X(\Q_{61})\to \Q_{61}$$
  in all affine
non-Weierstrass disks. 
By Lemma~\ref{L:X0N+model}, the local heights $h_{\ell}(\mathrm{A}(x))$ are trivial for
  $\ell\ne {61}$, so $\Upsilon =\{0\}$ and all rational points are
  zeroes of $\rho$.
We find that $\rho$ indeed vanishes on the known rational points, and that these are simple
zeroes of $\rho$. 

In addition, $\rho$ vanishes to multiplicity~1 on~82 points in $X(\Q_{61})$ that do not appear 
to be rational. 
As described in~\S\ref{subsec:mws}, these yield cosets of $61^2J(\Q)$, and our
implementation of the Mordell--Weil sieve shows that the image of these cosets does not intersect the
image of $X(\F_{229})$ inside $J(\F_{229})/61^2J(\F_{229})$.
Hence these additional zeroes do not come from a rational point.

Recall that there are no $\Q_{61}$-rational points at infinity, so
it only remains to show that there are no rational points in the Weierstrass disk.
To this end, we show that for 
$$S=\{41,83,641,1697,4057,10853\},$$
the image of the
reduction of this disk does not intersect $$\mathrm{im}(\beta_{S,\,2\cdot 61})\subset
\prod_{v\in S} \faktor{J(\F_v)}{MJ(\F_v)}\,,$$ 
where $M = 2\#J(\F_{61})$
and $\beta_{2,61}\colon \prod_{v\in S}X(\F_v)\to \prod_{v\in S}
J(\F_v)/MJ(\F_v)$ is induced by the Abel-Jacobi map with respect to $b$ and
the canonical surjections.

This completes the proof that $\#X(\Q)=6$. According to Galbraith~\cite{Gal96}, these points are all cusps or CM-points.
\end{Example}

\begin{Example}\label{E:167}
We were able to prove that the curve 
  \[
  X_0^+(167) \colon y^2 = x^6 - 4x^5 + 2x^4 - 2x^3 - 3x^2 + 2x - 3 
\]
only contains the four obvious rational points
$
  \{ (- 1, \pm1), \infty_{\pm}\}\,;
$
  these are all cusps or CM by Galbraith~\cite{Gal96}.
In our computation, we use our quadratic Chabauty algorithm for $p=7$ and the
  Mordell--Weil sieve, following the same strategy as in Example~\ref{E:107}.
  The verification that the additional solutions of the resulting $p$-adic
  functions are not rational was the most challenging Mordell--Weil sieve computation we
  encountered in our work; it required the auxiliary integer $5\cdot 11\cdot 19$ and the
  set of good primes
 $$S =\{3,5,19,29,31,67,263,281,283,769,1151,2377,3847,4957,67217\}.$$
\end{Example}

\begin{Example}\label{E:191}
  A model for $ X_0^+(191)$ is given by 
$$ y^2 =  x^6 + 2x^4 + 2x^3 + 5x^2 - 6x + 1\,.$$
We use quadratic Chabauty for $p=31$ together with the Mordell--Weil sieve exactly as above
to show that $X_0^+(191)(\Q) = \{  (0,\pm 1),  (2, \pm 11), \infty_{\pm}\}.$
  Galbraith (see~\cite[Table~1]{Gal99}) has shown that $(2, -11)$ is exceptional, with corresponding $j$-invariant
  \begin{align*}
    j = &2891249511562231668955764266428063102082570956800000 \\&\pm
  64074939271375546714155254091066566840131584000\sqrt{ 61 \cdot 229 \cdot
  145757}\,.
  \end{align*}
\end{Example}
%%%%%%%%%%%%%%%%%%%%%%%%%%%%%%%%%%%%%%%%%%%%%%%%%%%%%%%%%%%%%%%%%%%%%%%%
\subsubsection{Genus 3}\label{subsec:x0n+g3}

We apply our algorithm to show that the rational points on the curves $X_0^+(N)$ for $N$
as in~\eqref{eq:X0N+g3} are precisely the
ones already found by Galbraith.
All curves in our list are non-hyperelliptic and they have the convenient feature that they have sufficiently many
rational points, so no heights on divisors need to be computed.
We always find two independent cycles in $\ker(\NS(J)\longrightarrow \NS(X))$, and, as expected, the common zero set of the
corresponding functions consists precisely of the rational points found by Galbraith.

\begin{theorem}\label{T:X0N+g3}
  Let $N$ be a prime such that $X_0^+(N)$ has genus 3. Then the rational points on
  $X_0^+(N)$ are as below. In particular, all rational points are either cusps or
  CM-points, with discriminant $\Delta$.
\end{theorem}

\begin{Example}\label{E:97}
  A model for $ X_0^+(97)$ is given by 
\[
zx^3 + (-y^2 + zy)x^2 + (-y^3 - zy^2 - z^3)x + (zy^3 + z^2y^2) = 0.
\]
Using our algorithm for $p=5$, we find that
the rational points are as follows:

\tiny{\begin{tabular}{|c|c|c|c|c|c|c|c|c|c|c|}
\hline
$\Delta$ & cusp & $-3$ & $-4$ & $-8$ & $-11$ & $-12$ & $-16$ & $-27$ & $-43$ & $-163$\\
\hline
Point &  $(1: 0: 0)$ & $(-2: 1:1)$ & $(-1: 0:1)$ & $(0: 0: 1)$ & $(0: 1: 0)$ & $(0: -1:
  1)$  & $(1: 0: 1)$  & $(1: 1: 1)$ & $(-1: 1: 0)$ & $(5: 3: 2)$ \\
\hline
\end{tabular}}

\normalsize

\end{Example}

\begin{Example}\label{E:109}
  A model for $ X_0^+(109)$ is given by 
\[
zx^3 + (zy + z^2)x^2 + (-y^3 - zy^2 - z^3)x + (-zy^3 - 3z^2y^2 - 2z^3y) = 0.
\]
Using our algorithm for $p=29$, we find that
the rational points are as follows:

\tiny{\begin{tabular}{|c|c|c|c|c|c|c|c|c|c|c|}
\hline
$\Delta$ & cusp &$-3$ &  $-4$ & $-7$ & $-12$ &  $-16$ & $-27$ & $-28$ & $-43$ \\
\hline
Point & $(1: 0: 0)$ & $(-2: 1: 2)$ & $(0: -2: 1)$  & $(0: -1: 1)$  & $(0: 1: 0)$  & $(0:
  0: 1)$ & $(-1: -1: 1)$ & $(-2: 1: 1)$ & $(1: -1: 1)$ \\
\hline
\end{tabular}}
\normalsize
\end{Example}

\begin{Example}\label{E:113}
  A model for $ X_0^+(113)$ is given by 
\[
zx^3 + (-y^2 - z^2)x^2 + (y^3 + z^3)x + (-2z^2y^2 + z^3y) = 0.
\]
Using our algorithm for $p=17$, we find that
the rational points are as follows:

\tiny{\begin{tabular}{|c|c|c|c|c|c|c|c|c|c|}
\hline
$\Delta$ & cusp & $-4$ &  $-7$ & $-8$ & $-11$ & $-16$ &  $-28$ &  $-163$   \\
\hline
Point & $(1: 0: 0)$ & $(2: 2: 1)$  & $(0: 1: 0)$  & $(1: 1: 1)$  & $(1: 1: 0)$  & $(0: 0: 1)$  & $(0: 1: 2)$  & $(5: 3: 1)$ \\
\hline
\end{tabular}}
\normalsize
\end{Example}

\begin{Example}\label{E:127}
  A model for $ X_0^+(127)$ is given by 
\[
zx^3 + (-y^2 - 3z^2)x^2 + (y^3 - z^2y + 4z^3)x + (2zy^3 - 3z^2y^2 + 3z^3y - 2z^4) = 0.
\]
Using our algorithm for $p=11$, we find that
the rational points are as follows:

\tiny{\begin{tabular}{|c|c|c|c|c|c|c|c|c|}
\hline
$\Delta$ & cusp &   $-3$ &  $-7$ &  $-12$ &  $-27$ &  $-28$ &  $-43$ &  $-67$  \\
\hline
Point & $(1: 0: 0)$  & $(5: 3: 2)$  & $(2: 1: 1)$  & $(1: 1: 0)$  & $(1: 0: 1)$  & $(0: 1: 1)$  & $(0: 1: 0)$  & $(4: 2: 1)$ \\
\hline
\end{tabular}}
\normalsize
\end{Example}

\begin{Example}\label{E:139}
  A model for $ X_0^+(139)$ is given by 
\[
zx^3 + (-y^2 + zy)x^2 + (-y^3 - 2zy^2 - 3z^2y - z^3)x + (y^4 + zy^3 + z^2y^2 + z^3y) = 0.
\]
Using our algorithm for $p=19$, we find that
the rational points are as follows:

\tiny{
\begin{tabular}{|c|c|c|c|c|c|c|c|c|}
\hline
$\Delta$ & cusp&  $-3$ &  $-8$ &  $-12$ &  $-19$ &  $-27$ &  $-43$  \\
\hline
Point & $(1: 0: 0)$  & $(4: -3: 1)$ & $(0: 0: 1)$  & $(0: -1: 1)$  & $(1: -1: 1)$  & $(1: 0: 1)$  & $(-1: 0: 1)$ \\
\hline
\end{tabular}}
\normalsize
\end{Example}

\begin{Example}\label{E:149}
  A model for $ X_0^+(149)$ is given by 
\[
zx^3 - y^2x^2 + (y^3 + zy^2 - 2z^2y - z^3)x + (-y^4 + zy^3 + z^2y^2 - z^3y) = 0.
\]
Using our algorithm for $p=11$, we find that
the rational points are as follows:

\tiny{
\begin{tabular}{|c|c|c|c|c|c|c|c|}
\hline
$\Delta$ & cusp & $-4$ &  $-7$ &  $-16$ &  $-19$ & $-28$ &  $-67$  \\
\hline
Point & $(1: 0: 0)$  & $(-1: 0: 1)$  & $(0: 1: 1)$  & $(1: 0: 1)$ & $(0: 0: 1)$  & $(0: -1: 1)$  & $(2: 2: 1)$ \\
\hline
\end{tabular}}
\end{Example}

\begin{Example}\label{E:151}
  A model for $ X_0^+(151)$ is given by 
\[
zx^3 + (-2zy + z^2)x^2 + (-y^3 + 2zy^2)x + (-zy^3 + 3z^2y^2 - z^3y - 2z^4) = 0.
\]
Using our algorithm for $p=19$, we find that
the rational points are as follows:

\tiny{
\begin{tabular}{|c|c|c|c|c|c|c|c|c|}
\hline
$\Delta$ & cusp&  $-3$ & $-7$ & $-12$ &  $-27$ & $-28$ &  $-67$ &  $-163$ \\
\hline
Point & $(1: 0: 0)$ & $(-2: -2: 1)$  & $(0: 1: 0)$  & $(0: 2: 1)$  & $(1: 1: 1)$  & $(2: 3: 2)$ & $(1: 0: 1)$ & $(3: 2: 1)$ \\
\hline
\end{tabular}}
\end{Example}

\begin{Example}\label{E:179}
  A model for $ X_0^+(179)$ is given by 
\[
zx^3 + (-2zy - z^2)x^2 + (-y^3 - zy^2 - 2z^2y - z^3)x + (-zy^3 + z^3y) = 0.
\]
Using our algorithm for $p=17$, we find that
the rational points are as follows:

\tiny{\begin{tabular}{|c|c|c|c|c|c|c|}
\hline
$\Delta$ &  cusp&  $-7$ &  $-8$ & $-11$ &  $-28$ &  $-163$ \\
\hline
Point & $(1: 0: 0)$  & $(0: -1: 1)$ & $(0: 1: 0)$  & $(0: 0: 1)$  & $(0: 1: 1)$ & $(-2:  2: 1)$\\
\hline
\end{tabular}}
\end{Example}

\begin{Example}\label{E:239}
  A model for $ X_0^+(239)$ is given by 
\[
zx^3 + (-y^2 + zy + z^2)x^2 + (-y^3 - zy^2 - z^2y)x + (y^4 + 3zy^3 + 2z^2y^2 + z^3y) = 0.
\]
Using our algorithm for $p=13$, we find that
the rational points are as follows:

\tiny{
\begin{tabular}{|c|c|c|c|c|c|}
\hline
$\Delta$ & cusp &  $-7$ &  $-19$ &  $-28$ &  $-43$ \\
\hline
Point & $(1: 0: 0)$ & $(-1: 0: 1)$  & $(0: 0: 1)$  & $(1: -2: 1)$  & $(1: -1: 1)$ \\
\hline
\end{tabular}}
\end{Example}

%%%%%%%%%%%%%%%%%%%%%%%%%%%%%%%%%%%%%%%%%%%%%%%%%%%%%%%%%%%%%%%%%%%%%%%%
\subsection{Genus 2 curves in Mazur's Program B}\label{subsec:dzb}
In this section, we determine the rational points on two genus 2 curves that were
communicated to us by David Zureick-Brown. They arise in the
work of Rouse, Sutherland, and Zureick-Brown \cite{RSZB} on Mazur's
Program B as modular curves $X_H = X(25)/H$, where $\Gamma(25)\subset H\subset
\GL_2(\Z_5)$. Both curves have the following properties:
\begin{itemize}
  \item They each have 
two rational points of exponential height at
  most~1000, good reduction away from 5, and potentially good reduction at 5.
  \item 
Their Jacobians have real multiplication, no rational torsion and Mordell--Weil rank~2; they are both
    absolutely simple.
\item 
The Galois action on the 2-torsion field is $A_5$, which is too large for an elliptic curve
    Chabauty computation.
\end{itemize}
We prove that $\#X_H(\Q)=2$ for each curve $X_H$ using quadratic Chabauty and the Mordell--Weil
sieve, similar to the computation of $X^+_0(107)(\Q)$ described in detail in
Example~\ref{E:107}.
\begin{Example}\label{E:DZB1}
  A suitable affine model of the  curve $X_{11}$ is given by
  $$X_{11}\colon y^2 = -35x^6 + 310x^5 - 675x^4 + 750x^3 - 450x^2 + 140x - 15.$$
  As in Example~\ref{E:107}, we found the rather large prime $p=61$ to be the most convenient one for our
  computations.  
  We determine the height pairing on the Jacobian using divisors as
  in~\S\ref{subsec:coleman-gross}. The
  quadratic Chabauty function $\rho$ has 62 solutions in addition to the
  rational ones. Applying the
  Mordell--Weil sieve with the primes 7, 29, 257 and 3457, we show that these are in fact not rational; to prove 
  non-existence of rational points in the unique Weierstrass disk, we sieve with the
  primes 31, 61 and 191.
  This shows that $X_{11}(\Q) = \{(1, \pm 5)\}$. 
\end{Example}

\begin{Example}\label{E:DZB2}
  We use the model
\[ 
  X_{15} \colon y^2 =  5x^6 - 50x^4 - 150x^3 + 25x^2 + 90x + 25 
\]
  with small rational points $(0,\pm 5)$.
  Again we run quadratic Chabauty for a fairly large prime, namely $p=71$, resulting in~78
  additional
  zeroes in $X(\Q_{71})$ that we show to be non-rational by sieving with the primes
  $7,43,83,101,$ and  $1399$. There is an additional final sieving to show there are no
  rational points in the Weierstrass disk. We conclude that
    $X_{15}(\Q)=\{(0,\pm 5)\}$.
\end{Example}

%%%%%%%%%%%%%%%%%%%%%%%%%%%%%%%%%%%%%%%%%%%%%%%%%%%%%%%%%%%%%%%%%%%%%%%%
\subsection{Two curves with nontrivial local heights away from
$p$}\label{subsec:nontriv_away} 
We compute the rational points on two genus~2 curves $C_{188}$ and $C_{161}$ considered
in~\cite{FLSSSW}. In both cases, the Jacobian of $C_N$ is an optimal quotient of $J_0(N)$,
so it has real multiplication and Picard number~2. The Mordell--Weil ranks are both~2 as
well, and the rational torsion subgroup is trivial. In~\cite{FLSSSW} empirical evidence was
presented that the full conjecture of Birch and Swinnerton-Dyer holds for both Jacobians.
The curves themselves have good reduction away from $N$. 

So far, all curves whose rational points were computed via quadratic Chabauty had trivial
contributions away from $p$, except for the bielliptic examples in~\cite{BD18,BD19}.
However, for those examples it was possible to find the local contributions away from $p$ by relating them to local heights on
the elliptic quotients. In the examples presented here, we compute these contributions
using Theorem~\ref{thm:BeD}. As discussed in~\S\ref{subsec:betts_dogra}, we do not have a
general algorithm for the action induced by an endomorphism on \'etale cohomology. Nevertheless,
we show below that we can sometimes derive sufficient information from
Theorem~\ref{thm:BeD} to pin down the local contributions precisely, by computing the
local heights at $p=3$ for the known rational points
and by exploiting the bilinearity of the global height pairing.

We include these examples to illustrate the practicality of our algorithms. However, we note that the
rational points on both curves can be computed by combining covering collections with
elliptic curve Chabauty. For $C_{188}$ this was pointed out to us by Nils Bruin, and for
$C_{161}$, this computation is due to Bars, Gonz\'{a}lez, and Xarles~\cite{BGX20}.

\begin{Example}\label{E:188}
We first consider the genus 2 curve 
\begin{equation}\label{eq:188}
  C_{188}  \colon  y^2 = x^5 - x^4 + x^3 + x^2 - 2x + 1\,. 
\end{equation}
Over $\Z_{47}$,  it has a regular
semistable model whose special fibre is a curve of genus~1 with a node, so $h_{47}$ is
trivial by Theorem~\ref{thm:BeD}. However, as we shall see, there are nontrivial
contributions to the local height at 2.

  The integral points on  $C_{188}$   over $\Q(\sqrt{-3})$ were computed
  in~\cite[Example~6.5]{BBBM}.
 In the present work, we show that 
\begin{equation}\label{eq:188pts}
  C_{188}(\Q) = \{(0, \pm 1), (1, \pm 1), (-1, \pm 1), (2, \pm 5), (4, \pm 29),
  \infty\}.
\end{equation}
 For our computations, we use the good ordinary prime $p=3$, the base point $b=(1,1)$, and
  a cycle $Z$ constructed from the Hecke operator $T_3$ as in~\eqref{eq:tr0}.

  Recall from Example~\ref{E:h2_188} that there is a regular semistable model over 
  $K=\Q_2(\sqrt[3]{2})$ and that the corresponding  metric graph $\Gamma_{\reg}$ is a line
  segment.
The two genus one vertices $w_0 $ and $w_1$ have pre-images 
  $$\mathcal{U}_0 \colonequals \{P\in C_{188}(\Q _2 ):\ord_2 (x(P))>0\}\;,\;\mathcal{U}_1
  \colonequals \{P\in
  C_{188}(\Q _2 ):\ord_2 (x(P))=0\},$$ respectively. The set $\mathcal{U}_2 \colonequals \{P\in
  C_{188}(\Q _2 ):\ord_2 (x(P))<0 \}$ maps to the midpoint $w_2$ of the line segment. 
   
  Since the function $j_{\Gamma }$ from Theorem~\ref{thm:BeD} is affine linear and vanishes at $w_1$,
  there
  is a constant $\kappa$ such that  for all $x\in C_{188}(\Q_2)$ we have
  \[
  h_2 (\mathrm{A}(P))=m(P)\cdot {\kappa},
\]
where 
$$m(P)=
\begin{cases}
$2$, &\;\textrm{when $x(P)$ is divisible by 2},\\
$0$, &\;\textrm{when $x(P)$ is a $2$-adic unit},\\
$1$, &\;\textrm{when $x(P)$ is  non-integral}.
\end{cases}$$ One could determine $\kappa$ by further computing the
trace of $Z$ acting on the cohomology of the two genus one curves in the special fibre
  of the regular model described in Example~\ref{E:h2_188}. In this example, we can determine $\kappa$ by computing
local heights at $p$, as there is a unique choice of $\kappa$ such that $$h_3 (\mathrm{A}
(P))+m(P)\cdot \kappa$$ satisfies the bilinearity properties of a global height. We can reduce
the determination of 
  $\kappa$ to linear algebra by computing $h_3(\mathrm{A}(P))$ and the values of a basis of the space of
  $\End_0(J)$-equivariant bilinear pairings for $3=g+1$ rational points $P\in X(\Q)$.
  We find $\kappa=\frac{4}{3}\log_{p} (2)$. 

To finish the computation of the rational points, we first solve for the zeroes of the quadratic Chabauty function $\rho$ on the affine
patch~\eqref{eq:188}. In order to deal with the Weierstrass disk at
  infinity, we move the point at infinity to $(0,0)$ and repeat
  the computation for the resulting affine patch. We then apply the trick described
in~\cite[\S5.5]{BDMTV19}, changing the base point and reducing the computation of the Frobenius structure
to the computation of Coleman integrals. 

We find that $\rho$ vanishes on the known rational points and that it
  vanishes on~13
additional $\Q_3$-points to precision $3^5$. Upon noticing that $J(\F_{43})\simeq
(\Z/54\Z)^2$, we show that
  the reductions of the corresponding cosets of $27J(\Q)$ do not meet the image of
  $C_{188}(\F_{43})$ in $J(\F_{43})/27J(\F_{43})$. This
suffices to prove~\eqref{eq:188pts}.

\end{Example}

\begin{Example}\label{E:161}
  The curve $C_{161}$ has an affine equation
\begin{align*}
   y^2 & =x^6+2x^4+6x^3+17x^2+18x+5  = (x^3-2x^2 +3x+5)(x^3+2x^2 +3x+1).
\end{align*} 
  As discussed in~\cite{BGX20}, this is in fact a model for the modular curve $X_0^*(161)
  = X_0(161)/\langle w_7, w_{23}\rangle $.
The curve has ten small rational points
  \begin{equation}\label{eq:161pts}
\left(\frac{1}{4},\pm \frac{209}{64}\right),(-1,\pm 1), (1,\pm 7),\left(\frac{1}{2},\pm \frac{35}{8}\right),
    \infty_{\pm}.
  \end{equation}
Anticipating the need to use the Mordell--Weil sieve, we choose the prime $p=29$ and the
  cycle $Z$ corresponding to the endomorphism $4T_{29}-\Tr(T_{29})I_4$. 

The bad primes are $7$ and $23$. At both of these primes, the stable model has special
  fibre a genus zero curve with two double points.  One can show this, for
  instance, using the program {\tt genus2reduction} due to Qing Liu, now
  contained in {\tt Pari/GP} or {\tt Sage}.  This (or {\tt Magma}'s
  {\tt RegularModel} package) also shows that the model over $\Z_{23}$
  defined by the given equation is regular. Indeed, the 23-adic valuation
  of the discriminant is 2; therefore both singular points $(2,0)$ and
  $(11,0)$ on the reduction modulo~23 define regular points on this model. 
  Hence the given equation defines a
  regular semistable model over $\Z_{23}$, and all of the $\Q _{23}$ points lie on a common
  irreducible component of a minimal regular model over $\Z_{23}$, so the height
  contribution at this prime is zero by Theorem~\ref{thm:BeD}.

  At 7, 
  the discriminant has valuation~4, so the model defined by the given equation is not regular. The singular
  points on the special fibre are $(1,0)$ and $(4,0)$. Blowing up once in 
  both of these yields a semistable regular model whose special fiber 
consists of two genus~0 curves $w_1$ and $w_2$ that do not intersect and another genus~0
  curve $w_0$ which reduces to the smooth locus of the stable model and which
  intersects $w_1$ and $w_2$ transversely in two points each, $e_1 $ and $e_2 $ and $e_3 $ and $e_4$ respectively. 
  This information can also be obtained from {\tt genus2reduction} or {\tt
  RegularModel}.

\begin{figure}[htbp!]
\begin{center}
\begin{tikzpicture}
%%edges
\draw[thick] (1,0) circle (1);
\draw[thick] (3,0) circle (1);
 \node[] at (-0.4,0) {\footnotesize $w_1$};
 \node[] at (1.6,0) {\footnotesize $w_0$};
 \node[] at (4.4,0) {\footnotesize $w_2$};
 \node[] at (1.0,-1.3) {\footnotesize $e_1$};
 \node[] at (1.0,1.2) {\footnotesize $e_2$};
 \node[] at (3.0,-1.3) {\footnotesize $e_3$};
 \node[] at (3.0,1.2) {\footnotesize $e_4$};
%%vertices
  \coordinate (a) at (0,0);
\draw[fill=black] (a) circle (3pt);
\coordinate (b) at (2,0);
\draw[fill=black] (b) circle (3pt);
\coordinate (c) at (4,0);
\draw[fill=black] (c) circle (3pt);
%%genus of vertices
\end{tikzpicture}
\end{center}
  \caption{Dual graph of the minimal regular model of $C_{161}$ at $\ell=7$. }\label{fig1}
\end{figure}
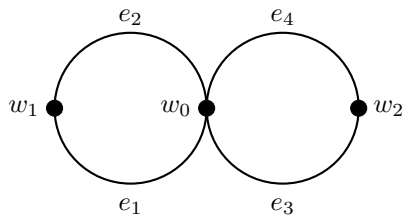

The corresponding dual graph is shown in Figure \ref{fig1}. 
We choose an orientation by designating $w_0$  as the source of $e_1 $ and $e_3 $ and as the target of $e_2 $ and $e_4$.
    The points $\left(\frac{1}{4},\pm
  \frac{209}{64}\right),(-1,\pm 1), \infty_{\pm}$ listed
  in~\eqref{eq:161pts} reduce to the component $w_0$. The points $(1,\pm 7)$
  reduce to $w_1$
  and the points $\left(\frac{1}{2},\pm \frac{35}{8}\right)$ reduce to
  $w_2$. 
We may again use Theorem \ref{thm:BeD} to determine the possible values of $h_7 (P)$,
  without computing the action of our chosen correspondence on $\HH^1 (\Gamma )$. Note that
  in this case, the homology $\HH_1 (\Gamma )$ is generated by $\gamma
  _1 =e_2 + e_1 $ and $\gamma _2 =e_3 +e_4$ respectively.
  Since $Z$ is trace zero on $\HH_1
  (\Gamma )$, 
  with respect to this basis, the corresponding endomorphism must be of the form
$\left(
\begin{array}{cc}
a & b \\
c & -a \\
\end{array}
\right).$  Then, by Theorem \ref{thm:BeD}, the measure $\mu _Z $ is simply given by $\frac{a}{2}(\gamma _1
  -\gamma _2 )$, since both edges have length 2. The image of $X(\Q _7 )$ in $\Gamma $ consists
  of 
  the three vertices $w_0, w_1,$ and $w_2$
  and hence if we take the
  basepoint $(\frac{1}{4},\frac{209}{64})$ reducing to $w_0$, 
  the values of $j_{\Gamma }$ are simply $a,0,-a$. 
  We solve for $a$ using a 29-adic computation similar to the previous example, and we find that $a=-4$.
Finally, we apply the Mordell--Weil sieve with $M=4\cdot 29^3$ and primes $199,373,463$ to
show that the only 29-adic points in the zero set of $\rho$ modulo $29^3$ are the rational
points listed in~\eqref{eq:161pts}. This proves that these are indeed the only rational
  points on $C_{161}$.
  \end{Example}

\subsection{The nonsplit Cartan modular curve $X_{\ns} ^+ (17)$}\label{sec:ns17}
The modular curve $$X \colonequals  X_{\rm ns} ^+ (17)$$ attached to the normaliser of the
non-split Cartan subgroup of level $17$ has genus $6$. By~\cite[\S5.3]{DLF19}, the rank of
$J_{\rm ns}^+(17)(\Q)$ is also~6. The set of rational points $X(\Q )$
can be determined \emph{without} computing local heights at the bad prime $17$, even
though these contribute nontrivially when determining $X(\Q_p)_2$, by choosing the
correspondence $Z$ carefully. 

\par The curve $X$ has a semistable model $\mathcal{X}$ over $W(\overline{\F}_{17})[\varpi]$ with $\varpi = (1+\zeta_{17})^{1/9}$, where $W(\overline{\F}_{17})$ is the ring of Witt vectors of $\overline{\F}_{17}$, described by Edixhoven--Parent \cite{EP21}. Its special fibre has two irreducible components 
\[\begin{array}{llll}
C_1 &:& y^2 = x(x^9 +a),  & \qquad a \in \overline{\F}_{17}^{\times}, \\
C_2 &:& z^2 = w(w^3 +b), & \qquad b \in \overline{\F}_{17}^{\times},
\end{array}\]
which have genus $4$ and $1$, respectively. They are smooth and intersect transversely in two points, so that the Jacobian has toric rank $1$. The inertia subgroup $I \subset G_{\Q_{17}}$ acts via automorphisms on the special fibre of this model, and the stabiliser of the set of irreducible components is contained in $\mu_{18}(\overline{\F}_{17^2}) = \langle \zeta \rangle$, where the root of unity $\zeta^{18}=1$ acts on the components by
\[\begin{array}{lllll}
\zeta &:& (x,y) & \longmapsto & (\zeta^4x,\zeta^2y), \\
\zeta &:& (z,w) & \longmapsto & (\zeta^{12}z,\zeta^{6}w).
\end{array}\]
The resulting operator $[\zeta]$ on the cohomology of these curves has characteristic polynomial 
\[
\begin{array}{lll}
\det\left(1 - t [\zeta] : \HH^1(C_1,\Q_p) \right) &=& (t^2+t+1)(t^6+t^3+1), \\
\det\left(1 - t [\zeta] : \HH^1(C_2,\Q_p) \right) &=& (t^2+t+1).
\end{array}
\]
Since the Hecke action on the cohomology of $X$ is defined over $\Q$, it must commute with
the action of inertia, and therefore the irreducible Hecke modules of the Jacobian up to
isogeny must be contained in the submodules coming from the toric part (dimension $1$) and
the parts where the operator $[\zeta]$ is of order $3$ (dimension $2$) and of order $9$
(dimension $3$). By the work of Chen and Edixhoven--de Smit \cite{Che00,EdS00} the
Jacobian of $X$ admits an isogeny to the new part of the Jacobian of $X_0^+(17^2)$
equivariant for the anemic Hecke algebra. The new part of the Jacobian of $X_0^+(17^2)$
decomposes up to isogeny into irreducible factors $M_1 \times M_2 \times M_3$ of
dimensions 
$1,2,3$ respectively, where $M_1,M_2,M_3$ are killed by the Hecke operators 
\[
\begin{array}{lll}
M_1 & : & (T_{2} + 1) = 0, \\
M_2 & : & (T^2_{2} +T_2 - 3) = 0,\\
M_3 & : & (T_{2}^3 -3T_{2} + 1) = 0. \\
\end{array}
\]
If we set $M = (T_{2} +1)(T^2_{2} +T_2 - 3)$, we find that $Z=M$ and $Z=2M^3 + 3M^2$ are
nontrivial trace zero correspondences that induce the zero endomorphisms on $H_1 (\Gamma ,\Q )$ and the cohomology of $C_2$,  so that Theorem \ref{thm:BeD} implies that $\mu _F = 0$, and hence the $17$-adic height vanishes:
$$h_{17}(\mathrm{A}_Z (x))=0, \qquad \textrm{for all}\; x\in X(\Q _{17}).$$
In fact, starting from any generator $T$ of the Hecke algebra (like $T=T_2$ above), one easily computes two linearly independent trace zero correspondences $Z \in \Z[T]$ that act trivially on the dual graph and the cohomology of $C_2$, which therefore likewise ensures the triviality of the associated $17$-adic height.

\par To put these observations into action, we choose $p=31$ and use the model of $X$
found by Mercuri and Schoof \cite[\S 6]{MS20} as an intersection of six quadrics in
$\PP^5$. Our strategy for finding a suitable singular plane curve model largely follows
\cite{AABCCKW}: To find a model
with small coefficients, we use the Magma function \texttt{Genus6PlaneCurveModel}, and then apply an automorphism of $\PP^2$ to ensure that there
are two rational points at infinity (this speeds up the computation of
the Hodge filtration (see~\cite[Section~4]{BDMTV19}) where one passes to a number field
over which the divisor at infinity splits completely). We obtain a singular plane curve
model $Q(x,y)=0$, where
\begin{align*}
5\cdot Q(x,y)&=5y^6 + \left(24x + 12\right)y^5 + \left(-495x^2 - 543x -153 \right)y^4 + \\
&\qquad \left(-1472x^3 - 2814x^2 - 1719x - 337 \right)y^3 +\\
&\qquad \left(-1686x^4 - 4875x^3 - 4761x^2 - 1902x - 263\right)y^2 + \\
&\qquad \left(-540x^5 - 2082x^4 - 2952x^3 - 1875x^2 - 535x  - 56 \right)y +\\
&\qquad  188x^6 + 534x^5 + 567x^4 + 284 x^3 + 70x^2 + 7x.
\end{align*}
The fact that $T_{31}$ generates the Hecke algebra can be checked from the
LMFDB page for newforms of weight two, level 289, trivial character and
Atkin--Lehner eigenvalue one \cite{LMFDB}. We compute two correspondences
$Z \in \Z[T_{31}]$ as above, and obtain a pair of power series in each
residue disk, whose common zeroes to precision $O(31^{20})$ correspond to
the rational points 
\[
\left\{ \left(-\frac{4}{9},\frac{1}{9}\right), \left(-\frac{2}{3},-\frac{1}{3}\right), \left(-\frac{1}{2},\frac{1}{2}\right), (0,0), (-1,0) ,\infty _1 ,\infty _2 \right\} \subset X(\Q)
\]
where $\infty _1 $ and $\infty _2 $ are the points $(1:-1:0) $ and
$(1:-\frac{1}{5}:0)$. Therefore, this must be the full set of rational
points $X(\Q)$. These were already found by Mercuri--Schoof~\cite[\S
6]{MS20}; they are all CM points and the corresponding discriminants are
$-3, -7, -11, -12, -27, -28, -163$. This proves Theorem~\ref{T:Xns17}.

\begin{Remark}
  It would be interesting to use the techniques of this paper to compute the rational
  points on $X_{\ns} ^+ (19)$. Mercuri and
  Schoof~\cite[\S 7]{MS20} found a model for this curve as well. Nevertheless, we were unable to find a plane affine equation for this
  curve and a prime $p$, satisfying Assumption~\ref{tuitman1}, such that it is feasible
  to carry out Algorithm~\ref{alg:qc}. Difficulties arose in computing a basis of $\HH
  ^1 _{\dR}(X_{\Q _p })$ due to the large degrees of the field extensions we encountered
  when applying the algorithms in~\cite[\S 3]{Tui17}.
\end{Remark}
%=========================================================================

\bibliographystyle{alpha}

\begin{thebibliography}{BDCKW18}

\bibitem[AAB{\etalchar{+}}21]{AABCCKW}
N.~Ad{\v{z}}aga, V.~Arul, L.~Beneish, M.~Chen, S.~Chidambaram, T.~Keller, and
  B.~Wen.
\newblock Quadratic {C}habauty for {A}tkin-{L}ehner quotients of modular curves
  of prime level and genus 4, 5, 6.
\newblock {\em ArXiv preprint}, arXiv:2105.04811, 2021.



\bibitem[ACKP]{ACKP}
N.~Ad{\v{z}}aga, S.~Chidambaram, T.~Keller, and O.~Padurariu.
  \newblock {Rational points on hyperelliptic Atkin-Lehner quotients of modular curves and their coverings},
\newblock {\em Res. Number Theory}, 8:Art. 87, 2022.

\bibitem[AM]{AM}
 V.~Arul and J.\thinspace{}S. M\"uller
    \newblock Rational points on $X_0^+(125)$
\newblock {\em Expo. Math.}, to appear.

\bibitem[BB12]{bb:heights}
J.\thinspace{}S. Balakrishnan and A.~Besser.
\newblock Computing local {$p$}-adic height pairings on hyperelliptic curves.
\newblock {\em IMRN}, 2012(11):2405--2444, 2012.

\bibitem[BB15]{BB15}
J.\thinspace{}S. Balakrishnan and A.~Besser.
\newblock Coleman-{G}ross height pairings and the {$p$}-adic sigma function.
\newblock {\em J. Reine Angew. Math.}, 698:89--104, 2015.

\bibitem[BB21]{bb:heightserrata}
J.\thinspace{}S. Balakrishnan and A.~Besser.
\newblock Errata for ``{C}omputing local {$p$}-adic height pairings on
  hyperelliptic curves".
\newblock \url{http://math.bu.edu/people/jbala/cg_heights_errata.pdf}, 2021.

\bibitem[BBBM21]{BBBM}
J.\thinspace{}S. Balakrishnan, A.~Besser, F.~Bianchi, and J.\thinspace{}S.
  M\"uller.
\newblock Explicit quadratic {C}habauty over number fields.
\newblock {\em Israel J. Math}, 243:185--232, 2021.

\bibitem[BBM16]{BBM16}
J.\thinspace{}S. Balakrishnan, A.~Besser, and J.\thinspace{}S. M\"uller.
\newblock Quadratic {C}habauty: {$p$}-adic heights and integral points on
  hyperelliptic curves.
\newblock {\em J. Reine Angew. Math.}, 720:51--79, 2016.

\bibitem[BBM17]{BBM17}
J.\thinspace{}S. Balakrishnan, A.~Besser, and J.\thinspace{}S. M\"{u}ller.
\newblock Computing integral points on hyperelliptic curves using quadratic
  {C}habauty.
\newblock {\em Math. Comp.}, 86(305):1403--1434, 2017.

\bibitem[BBB{\etalchar{+}}21]{BBBLMTV19}
J.\thinspace{}S. Balakrishnan, A.J. Best, F.~Bianchi, B.~Lawrence,
  J.\thinspace{}S. M\"uller, N.~Triantafillou, and J.~Vonk.
\newblock Two recent {$p$}-adic approaches towards the (effective) {M}ordell
  conjecture.
\newblock In {\em Regulators {IV}: {A}n international conference on arithmetic
  {L}-functions and differential geometric methods}, volume 338 of {\em Progr.
  Math.}, pages 31--74. Birkh\"auser Boston, Boston, MA, 2021.
  
  \bibitem[BD18]{BD18}
J.\thinspace{}S. Balakrishnan and N.~Dogra.
\newblock Quadratic {C}habauty and rational points {I}: $p$-adic heights.
\newblock {\em Duke Math. J.}, 167(11):1981--2038, 2018.

\bibitem[BD21]{BD19}
J.\thinspace{}S. Balakrishnan and N.~Dogra.
\newblock Quadratic {C}habauty and rational points {II}: {G}eneralised height
  functions on {S}elmer varieties.
\newblock {\em Int. Math. Res. Not. IMRN}, (15):11923--12008, 2021.

\bibitem[BDM{\etalchar{+}}]{QCCode}
J.\thinspace{}S. Balakrishnan, N.~Dogra, J.\thinspace{}S. M\"uller, J.~Tuitman,
  and J.~Vonk.
\newblock {QCMod} (Magma code).
\newblock \url{https://github.com/steffenmueller/QCMod}.

\bibitem[BDM{\etalchar{+}}19]{BDMTV19}
J.\thinspace{}S. Balakrishnan, N.~Dogra, J.\thinspace{}S. M\"uller, J.~Tuitman,
  and J.~Vonk.
\newblock Explicit {C}habauty--{K}im for the split {C}artan modular curve of
  level $13$.
\newblock {\em Annals of Math.}, 189(3), 2019.

\bibitem[BDCKW18]{BDCKW}
J.\thinspace{}S. Balakrishnan, I.~Dan-Cohen, M.~Kim, and S.~Wewers.
\newblock A non-abelian conjecture of {T}ate-{S}hafarevich type for hyperbolic
  curves.
\newblock {\em Math. Ann.}, 372(1-2):369--428, 2018.

\bibitem[BKK11]{BKK10}
J.\thinspace{}S. Balakrishnan, K.\thinspace{}S. Kedlaya, and M.~Kim.
\newblock Appendix and erratum to ``{M}assey products for elliptic curves f
  rank 1''.
\newblock {\em J. Amer. Math. Soc.}, 24(1):281--291, 2011.

\bibitem[BT20]{BT19}
J.\thinspace{}S. Balakrishnan and J.~Tuitman.
\newblock Explicit {C}oleman integration for curves.
\newblock {\em Math. Comp.}, 89(326):2965--2984, 2020.

\bibitem[BC14]{banwait-cremona}
B.\thinspace{}S. Banwait and J.\thinspace{}E. Cremona.
\newblock Tetrahedral elliptic curves and the local-global principle for
  isogenies.
\newblock {\em Algebra \& Number Theory}, 8(5):1201--1229, 2014.

\bibitem[BGX21]{BGX20}
F.~Bars, J.~Gonz\'{a}lez, and X.~Xarles.
\newblock Hyperelliptic parametrizations of {$\Bbb{Q}$} curves.
\newblock {\em Ramanujan J.}, 56(1):103--120, 2021.

\bibitem[Bes04]{Bes04}
A.~Besser.
\newblock The $p$-adic height pairings of {C}oleman--{G}ross and of
  {N}ekov\'a\v{r}.
\newblock In {\em Number Theory}, volume~36 of {\em CRM Proc. Lect. Notes},
  pages 13--25. Amer. Math. Soc., 2004.


\bibitem[BMS]{BMS21}
  A.~Besser, J.\thinspace{}S. M\"uller and P.~Srinivasan.
\newblock $p$-adic adelic metrics and Quadratic Chabauty~I
\newblock {\em Arxiv preprint}, arXiv:2112.03873, 2021.
  
  \bibitem[BO83]{BO83}
P.~Berthelot and A.~Ogus.
\newblock F-isocrystals and de {Rham} cohomology {I}.
\newblock {\em Invent. Math.}, 72:159--199, 1983.
  

\bibitem[BD19]{BeD19}
A.~Betts and N.~Dogra.
\newblock The local theory of unipotent {K}ummer maps and refined {S}elmer
  schemes.
\newblock {\em ArXiv preprint}, arXiv:1909.05734v2, 2019.


\bibitem[vBHM20]{van2020explicit}
R.~van Bommel, D.~Holmes, and J.\thinspace{}S. M{\"u}ller.
\newblock Explicit arithmetic intersection theory and computation of
  {N}{\'e}ron-{T}ate heights.
\newblock {\em Math. Comp.}, 89(321):395--410, 2020.

\bibitem[BCP97]{BCP97}
W.~Bosma, J.~Cannon, and C.~Playoust.
\newblock {The Magma algebra system I: The user language}.
\newblock {\em J. Symb. Comp}, 24(3-4):235--265, 1997.


\bibitem[Box21]{Box19}
J.~Box.
\newblock Quadratic points on modular curves with infinite {M}ordell-{W}eil
  group.
\newblock {\em Math. Comp.}, 90(327):321--343, 2021.

\bibitem[BS10]{Bruin-Stoll:MWSieve}
N.~Bruin and M.~Stoll.
\newblock The {M}ordell-{W}eil sieve: proving non-existence of rational points
  on curves.
\newblock {\em LMS J. Comput. Math.}, 13:272--306, 2010.


\bibitem[CG89]{CG89}
R.\thinspace{}F. Coleman and B.\thinspace{}H. Gross.
\newblock {$p$}-adic heights on curves.
\newblock In {\em Algebraic number theory}, volume~17 of {\em Adv. Stud. Pure
  Math.}, pages 73--81. Academic Press, Boston, MA, 1989.

\bibitem[Che00]{Che00}
I.~Chen.
\newblock On relations between {J}acobians of certain modular curves.
\newblock {\em J. Algebra}, 231(1):414--448, 2000.


\bibitem[CR91]{CR91}
T.~Chinburg and R.~Rumely.
 \newblock Well-adjusted models for curves over {D}edekind rings
\newblock {\em Arithmetic algebraic geometry ({T}exel, 1989)}, Progr. Math., 89, 3--24, 1991.

\bibitem[CR93]{CR93}
T.~Chinburg and R.~Rumely.
\newblock The capacity pairing
\newblock {\em J. Reine Angew. Math.}, 434, 1993, 1--44.

\bibitem[Cla03]{clarkthesis}
P.\thinspace{}L. Clark.
\newblock {\em Rational points on {A}tkin-{L}ehner quotients of {S}himura
  curves}.
\newblock ProQuest LLC, Ann Arbor, MI, 2003.
\newblock Thesis (Ph.D.)--Harvard University.

\bibitem[CMSV19]{CMSV}
E.~Costa, N.~Mascot, J.~Sijsling, and J.~Voight.
\newblock Rigorous computation of the endomorphism ring of a {J}acobian.
\newblock {\em Math. Comp.}, 88(317):1303--1339, 2019.

\bibitem[DR73]{DR73}
P.~Deligne and M.~Rapoport.
\newblock Les schemas de modules de courbes elliptiques. 
\newblock {\em Modular Functions of one Variable II}, Proc. internat. Summer School, Univ. Antwerp 1972, Lect. Notes Math. 349, 143--316, 1973.

\bibitem[DF21]{DLF19}
N.~Dogra and S.~Le Fourn.
\newblock {Q}uadratic {C}habauty for modular curves and modular forms of rank
  one.
\newblock {\em Math. Ann.}, 380(1-2):393--448, 2021.


\bibitem[DRHS]{DRHS}
J. Duque-Rosero  and S. Hashimoto and P. Spelier.
\newblock Geometric Quadratic Chabauty and $p$-adic heights.
\newblock {\em Arxiv preprint}, arXiv:2207.10389, 2022.

\bibitem[EdS00]{EdS00}
B.~Edixhoven and B.~de~Smit.
\newblock Sur un r\'esultat d'{I}min {C}hen.
\newblock {\em Math. Res. Lett.}, (2--3):147--153, 2000.

\bibitem[EL21]{EL19}
B.~Edixhoven and G.~Lido.
\newblock Geometric quadratic {C}habauty.
\newblock {\em J. Inst. Math. Jussieu}, to appear,
  {\url{https://doi.org/10.1017/S1474748021000244}}, 2021.

\bibitem[EP21]{EP21}
B.~Edixhoven and P.~Parent.
\newblock Semistable reduction of modular curves associated with maximal
  subgroups in prime level.
\newblock {\em Doc. Math.}, 26:231--269, 2021.

\bibitem[FLS{\etalchar{+}}01]{FLSSSW}
E.\thinspace{}V. Flynn, G.~Lepr\'{e}vost, E.\thinspace{}F. Schaefer,
  W.\thinspace{}A. Stein, M.~Stoll, and J.\thinspace{}L. Wetherell.
\newblock Empirical evidence for the {B}irch and {S}winnerton-{D}yer
  conjectures for modular {J}acobians of genus 2 curves.
\newblock {\em Math. Comp.}, 70(236):1675--1697, 2001.


\bibitem[Gaj22]{Gaj22}
  S.~Gajovi\'{c}.
\newblock {\em Variations on the method of Chabauty and Coleman}
\newblock Thesis (Ph.D.)--University of Groningen, 2022.
\newblock
{\url{https://research.rug.nl/en/publications/variations-on-the-method-of-chabauty-and-coleman}}

\bibitem[Gal96]{Gal96}
S.\thinspace{}D. Galbraith.
\newblock Equations for modular curves.
\newblock {\em Oxford DPhil thesis}, 1996.

\bibitem[Gal99]{Gal99}
S.\thinspace{}D. Galbraith.
\newblock Rational points on {$X_0^+(p)$}.
\newblock {\em Experiment. Math.}, 8(4):311--318, 1999.

\bibitem[Gal02]{Gal02}
S.\thinspace{}D. Galbraith.
\newblock Rational points on {$X^+_0(N)$} and quadratic {$\Bbb Q$}-curves.
\newblock {\em J. Th\'eor. Nombres Bordeaux}, 14(1):205--219, 2002.

\bibitem[SGA7]{SGA7}
A.~Grothendieck
SGA 7, expos\'e IX. 
\newblock In {\em Lecture Notes in Mathematics 288}, pages 313-523. 
\newblock Springer-Verlag, New York, 1972.


\bibitem[Hol12]{Hol14}
D.~Holmes.
\newblock Computing {N}\'eron-{T}ate heights of points on hyperelliptic
  {J}acobians.
\newblock {\em J. Number Theory}, 132(6):1295--1305, 2012.

\bibitem[HZ02]{HZ02}
E.\thinspace{}W. Howe and H.\thinspace{}J. Zhu.
\newblock On the existence of absolutely simple abelian varieties of a given
  dimension over an arbitrary field.
\newblock {\em J. Number Theory}, 92(1):139--163, 2002.

\bibitem[Kim05]{Kim05}
M.~Kim.
\newblock The motivic fundamental group of $\mathbf{P}^1 \backslash \{
  0,1,\infty \}$ and the theorem of {S}iegel.
\newblock {\em Invent. Math.}, 161:629--656, 2005.

\bibitem[Kim09]{Kim09}
M.~Kim.
\newblock The unipotent {A}lbanese map and {S}elmer varieties for curves.
\newblock {\em Publ. RIMS}, 45:89--133, 2009.

\bibitem[Kim10]{Kim10b}
M.~Kim.
\newblock Massey products for elliptic curves of rank $1$.
\newblock {\em J. Amer. Math. Soc.}, 23(3):725--747, 2010.

\bibitem[KT08]{KT08}
M.~Kim and A.~Tamagawa.
\newblock The $l$-component of the unipotent {A}lbanese map.
\newblock {\em Math. Ann.}, 340(1):223--235, 2008.

\bibitem[Lig77]{Lig76}
G.~Ligozat.
\newblock Courbes modulaires de niveau $11$.
\newblock In {\em Modular Functions of One Variable V}, volume 601 of {\em
  Lecture Notes in Math.}, pages 149--237. Springer, Berlin, 1977.

\bibitem[LMFDB]{LMFDB}
The LMFDB Collaboration.
\newblock The {L}-functions and modular forms database, 2021.
\newblock \url{http://www.lmfdb.org}

\bibitem[Maz77]{maz77}
B.~Mazur.
\newblock Rational points on modular curves.
\newblock In {\em Modular functions of one variable, {V} ({P}roc. {S}econd
  {I}nternat. {C}onf., {U}niv. {B}onn, {B}onn, 1976)}, pages 107--148. Lecture
  Notes in Math., Vol. 601, 1977.

\bibitem[MS20]{MS20}
P.~Mercuri and R.~Schoof.
\newblock Modular forms invariant under non-split {C}artan subgroups.
\newblock {\em Math. Comp.}, 89(324):1969--1991, 2020.

\bibitem[M{\"u}l14]{Mul14}
J.\thinspace{}S M{\"u}ller.
\newblock Computing canonical heights using arithmetic intersection theory.
\newblock {\em Math. Comp.}, 83(285):311--336, 2014.

\bibitem[Nek93]{Nek93}
J.~Nekov{\'a}{\v{r}}.
\newblock On $p$-adic height pairings.
\newblock In {\em S\'eminaire de {T}h\'eorie des {N}ombres, {P}aris
  1990--1991}, pages 127--202. Birkh\"auser Boston, 1993.

\bibitem[PY07]{PY}
P.~Parent and A.~Yafaev.
\newblock Proving the triviality of rational points on {A}tkin-{L}ehner
  quotients of {S}himura curves.
\newblock {\em Math. Ann.}, 339(4):915--935, 2007.

%\bibitem[Ray90]{Ray90}
%M.~Raynaud.
%\newblock $p$-{G}roupes et r\'eduction semi-stable des courbes.
%\newblock In P.~Cartier, editor, {\em The Grothendieck Festschrift, vol. III},
%  volume~88 of {\em Progr. Math.}, pages 179--197. Birkh\"auser, 1990.

\bibitem[RSZB21]{RSZB}
J.~Rouse, A.\thinspace{}V. Sutherland, and D.~Zureick-Brown.
\newblock {$\ell$}-adic images of {G}alois for elliptic curves over
  {$\mathbb{Q}$} (and an appendix with J. Voight).
\newblock {\em Forum of Math. Sigma}, 10, E62, 2022.

\bibitem[RZB15]{RZB}
J.~Rouse and D.~Zureick-Brown.
\newblock Elliptic curves over {$\Bbb Q$} and 2-adic images of {G}alois.
\newblock {\em Res. Number Theory}, 1:Art. 12, 34, 2015.

\bibitem[Sch99]{Scharaschkin:Thesis}
V.~Scharaschkin.
\newblock {\em Local-global problems and the {B}rauer-{M}anin obstruction}.
\newblock ProQuest LLC, Ann Arbor, MI, 1999.
\newblock Thesis (Ph.D.)--University of Michigan.

\bibitem[Ser72]{Ser72}
J.-P. Serre.
\newblock Propri\'et\'es galoisiennes des points d'ordre fini des courbes
  elliptiques.
\newblock {\em Invent. Math.}, 15(4):259--331, 1972.

\bibitem[Sik17]{Sik17}
S.~Siksek.
\newblock Quadratic {C}habauty for modular curves.
\newblock {\em Arxiv preprint}, arXiv:1704.00473, 2017.

\bibitem[Tui16]{Tui16}
J.~Tuitman.
\newblock Counting points on curves using a map to {$\mathbf{P}^1$}.
\newblock {\em Math. Comp.}, 85(298):961--981, 2016.

\bibitem[Tui17]{Tui17}
J.~Tuitman.
\newblock Counting points on curves using a map to {$\bold{P}^1$}, {II}.
\newblock {\em Finite Fields Appl.}, 45:301--322, 2017.


\bibitem[Wal11]{Waldschmidt}
M.~Waldschmidt.
\newblock On the {$p$}-adic closure of a subgroup of rational points on an
  {A}belian variety.
\newblock {\em Afr. Mat.}, 22(1):79--89, 2011.

\bibitem[Xue09]{Xue09}
H.~Xue.
\newblock Minimal resolution of Atkin–Lehner quotients of $X_0(N)$. 
\newblock {\em J. Number Theory}, 129(9):2072--2092, 2009.

\bibitem[Zha93]{zhang}
S.~Zhang.
\newblock Admissible pairing on a curve.
\newblock {\em Invent. Math.}, 112(1):171--193, 1993.

\end{thebibliography}
\newcommand{\etalchar}[1]{$^{#1}$}

%=========================================================================
\end{document}